\documentclass[11pt,letterpaper]{amsart}
\usepackage[english]{babel}
\usepackage[utf8]{inputenc}
\usepackage{lmodern}
\usepackage[T1]{fontenc}
\usepackage{graphicx}
\usepackage{amsmath,amsthm,amssymb}
\usepackage{amsfonts}
\usepackage{verbatim}
\usepackage{amssymb}
\usepackage{enumitem}
\usepackage{stmaryrd}
\usepackage{esint}
\usepackage{soul}
\usepackage{floatrow}
\usepackage{subeqnarray}
\usepackage[all]{xy}
\usepackage{varioref}
\usepackage{tikz-cd}
\usepackage{hyperref}
\usepackage{pgfplots}

\newcommand{\proj}{\text{proj}}
\newcommand{\germ}{\text{germ}}
\newcommand{\Opp}{\text{Opp}}

\newcommand{\bdinf}{X^\infty}
\newcommand{\ch}{\text{ch}}
\newcommand{\chinf}{\text{ch}(X^\infty)}
\newcommand{\dinf}{\widetilde{d^\infty}}
\newcommand{\cat}{\text{CAT}}

\newcommand{\id}{\text{id}}

\newcommand{\Min}{\text{Min}}

\newcommand{\Aut}{\mathrm{Aut}}

\DeclareMathOperator{\supp}{\text{supp}}

\DeclareMathOperator{\iso}{Isom}
\newcommand{\bbR}{\mathbb{R}}
\newcommand{\bbP}{\mathbb{P}}

\newcommand{\bfH}{\mathbf{H}}

\newcommand{\fraka}{\mathfrak{a}}

\newcommand{\bbK}{\mathbb{K}}
\newcommand{\calF}{\mathcal{F}}
\newcommand{\calA}{\mathcal{A}}

\newcommand{\calD}{\mathcal{D}}
\newcommand{\calT}{\mathcal{T}}

\newcommand{\calO}{\mathcal{O}}

\newcommand{\calB}{\mathcal{B}}
\newcommand{\calU}{\mathcal{U}}

\newcommand{\SL}{\mathrm{SL}}

\newcommand{\stab}{\mathrm{Stab}}
\newcommand{\Id}{\mathrm{Id}}
\newcommand{\bd}{\partial_\infty}
\newcommand{\nui}{\check{\nu}}
\newcommand{\mui}{\check{\mu}}
\newcommand{\Bi}{\check{B}}

\newcommand{\prob}{\text{Prob}}

\newcommand{\Ch}{\mathrm{ch}}

\theoremstyle{plain}
\newtheorem{thm}{Theorem}[section]

\newtheorem{cor}[thm]{Corollary}
\newtheorem{lem}[thm]{Lemma}
\newtheorem{prop}[thm]{Proposition}
\theoremstyle{definition} 
\newtheorem{Def}[thm]{Definition}
\theoremstyle{remark}
\newtheorem{rem}[thm]{Remark}

\newtheorem{ex}[thm]{Example}

\newtheorem*{ackn}{Acknowledgements}

\newcommand{\newcomment}[4]{%
	\newcounter{#2counter}
	\expandafter\newcommand\csname #1\endcsname[1]{%
		\refstepcounter{#2counter}%
		{\color{#4}(#3\arabic{#2counter})}\marginpar{\scriptsize\raggedright\textbf{\color{#4}(#2 \arabic{#2counter}):} ##1}%
}}
\newcomment{clb}{Corentin}{c}{teal}

\newcommand{\blue}[1]{{\color{blue}#1}}
\newcommand{\red}[1]{{\color{red}#1}}

\title[Dynamical boundaries of affine buildings]{Dynamical boundaries of affine buildings: $C^\ast$-simplicity and Poisson boundaries}
\author{Corina Ciobotaru}
\author{Corentin Le Bars}
\date{}

\begin{document}
	\begin{abstract}
		We investigate a class of groups acting on possibly exotic affine buildings $X$ and possessing good proximal properties. Such groups are termed \emph{of general type}, and their dynamics is analyzed through their flag limit sets in the space of chambers at infinity of $X$. For a group $G$ of general type, we prove $C^\ast$-simplicity by showing that its flag limit set $\Lambda_{\mathcal F}(G)$ is topologically free, minimal, and strongly proximal. When $\Lambda_{\mathcal F}(G)$ intersects all Schubert cells relative to a limit chamber, then $\Lambda_{\mathcal F}(G)$ is a mean proximal space, in the sense that it carries a unique proximal stationary measure for any admissible probability measure on the acting group. Lattices are established as examples of groups of general type, and their Poisson boundaries are identified. The arguments rely on constructing an equivariant barycenter map from triples of chambers in generic position to the affine building.
	\end{abstract}
	\maketitle
	\setcounter{tocdepth}{1}
	\tableofcontents
	\section{Introduction}
	In this paper, we study groups acting on affine buildings and on their boundaries. Our approach relies on proximal properties of certain isometries of affine buildings that exhibit a weak form of north--south dynamics on the chambers at infinity. 
	These isometries, called \emph{strongly regular hyperbolic} in \cite{caprace_ciobotaru15}, behave similarly to loxodromic (i.e.\ $\mathbb{R}$-regular) elements of semisimple  real Lie groups.
	
	We further develop a theory of barycenters for triples of chambers at infinity, analogous to the barycenters associated with triples of distinct points in hyperbolic spaces. These tools are applied in two directions: to derive operator-algebraic properties of the groups under consideration, and to analyze their boundary theory.
	
	The level of generality in which we work allows us to formulate these results uniformly for a class of groups that we call \emph{of general type}, a class that includes in particular all cocompact lattices in affine buildings.

	\subsection{Groups of general type and their flag limit sets}	
	In this paper, $G$ denotes a discrete countable group acting by automorphisms on a thick, locally finite affine building $X$. We denote by $\chinf$ the set of chambers of the spherical building at infinity $X^\infty$ of $X$. An element $g \in G$ is said to be \emph{strongly regular hyperbolic} if it is an axial isometry whose translation axes are regular and thus lie in a single apartment of $X$; see Section~\ref{section SRH}. Such an element admits a pair of attracting and repelling chambers $C^\pm \in \chinf$, and successive powers of $g$ contract ``most'' chambers at infinity toward the attracting chamber $C^+$; see Proposition~\ref{prop dyn srh}.
	
	The \emph{flag limit set} $\Lambda_{\mathcal F}(G)$ of $G$ is defined as the closure, in the cone topology on $\chinf$, of the set $\Lambda^+_{\mathcal F}(G)$ of attracting chambers of strongly regular hyperbolic elements of $G$. In the linear case, this object was introduced by Benoist \cite[\S3]{benoist97}. We say that $G$ is of \emph{general type} if it contains strongly regular hyperbolic elements and if there are no ``large'' closed proper $G$-invariant subsets of $\chinf$; see Section~\ref{section general type} for a precise definition. This notion is designed so that if $G$ is a Zariski-dense subgroup of a semisimple algebraic group defined over a field $\bbK$, and $G$ contains a  loxodromic element, then $G$ is of general type when acting on the associated symmetric space or Bruhat-Tits building and on its spherical building at infinity. We point out that if the field $\bbK$ is Archimedean, then $G$ automatically contains loxodromic elements \cite{abels_margulis_soifer}, but this is no longer true if $\bbK$ is non-Archimedean \cite{quint02}.

	The following result shows that the flag limit set is a natural dynamical object. 
	
	\begin{thm}\label{thm flag general type intro}
		Let $X$ be a locally finite affine building and let $G < \iso(X)$ be a group of general type. Then the limit set $\Lambda_\calF(G)$ is the smallest closed non-empty $G$-invariant subset of $\chinf$. It is $G$-minimal, perfect and uncountable. If $\Lambda_\calF(G) \neq \chinf$, then $\Lambda_\calF(G)$ has empty interior. 
	\end{thm}
	
	In the algebraic case, such a statement is reminiscent of results due to Benoist, see for instance \cite[\S2.8.3]{benoist_quint}. 
	
	\subsection{$C^\ast$-simplicity for groups acting on buildings}
	
	Let $G$ be a discrete group and let $\lambda$ be the left regular representation of $G$ on the Hilbert space $\ell^2(G)$. The reduced $C^\ast$-algebra $C^\ast_r(G)$ is the closure of $\operatorname{span}\{\lambda(g) \mid g \in G\}\subseteq B(\ell^2(G))$ for the norm-operator topology. We say that $G$ is \emph{$C^\ast$-simple} if the $C^\ast$-algebra $C^\ast_r(G)$ does not have closed two-sided non-trivial ideals. An alternative characterization is that every representation of $G$ that is weakly contained in the left regular representation is actually weakly equivalent to the left regular representation \cite{de-la-harpe07}. 
	The first main result of this paper is the following. 
	
	\begin{thm}[{$C^\ast$-simplicity}]\label{thm c-star simpl intro}
		Let $G < \iso(X)$ be a group of general type acting properly on a locally finite affine building $X$. Then $G$ is $C^\ast$-simple.
	\end{thm}
	
	This result was previously known for lattices in Bruhat--Tits buildings \cite{bekka_cowling_de-la-harpe94}, and was later extended by Poznansky to discrete linear groups with trivial amenable radical \cite{poznansky09}; see \cite[Corollary~6.10]{breuillard_kalantar_kennedy_ozawa17}. It is also known to hold for triangle groups, which act simply transitively on affine buildings of type $\tilde{A}_2$ \cite{robertson_steger96}.
	
	We provide a uniform proof that does not require the action to be free or transitive on chambers, and which moreover applies to the exotic buildings of types $\tilde{C}_2$ and $\tilde{G}_2$. In the same way that Zariski-dense subgroups of semisimple Lie groups form a much larger class than lattices, the class of groups of general type is significantly broader than that of cocompact lattices in affine buildings. By \cite{bader_caprace_lecureux19}, exotic $\tilde{A}_2$-lattices are non-linear, and in fact all exotic lattices (of types $\tilde{C}_2$ and $\tilde{G}_2$) are conjectured to be non-linear. In such cases, the methods of \cite{bekka_cowling_de-la-harpe94} and \cite{poznansky09} do not apply.		
	
	The proof of Theorem~\ref{thm c-star simpl intro} is based on the dynamical characterization of $C^\ast$-simplicity established by Breuillard, Kalantar, Kennedy, and Ozawa in \cite{breuillard_kalantar_kennedy_ozawa17}. In particular, we construct a $G$-boundary on which $G$ acts topologically freely. Recall that a compact metrizable space $Y$ is called a $\Gamma$-boundary for a topological group $\Gamma$ if the action $\Gamma \curvearrowright Y$ by homeomorphisms is both strongly proximal and minimal. These notions are reviewed in Section~\ref{section furst bd}. Theorem~\ref{thm c-star simpl intro} then follows directly from the result below:
	
	\begin{thm}\label{thm G-boundary intro}
		Under the same assumptions as in Theorem~\ref{thm c-star simpl intro}, the flag limit set $\Lambda_\mathcal{F}(G) \subseteq \chinf$, equipped with the cone topology inherited from $\chinf$, is a $G$-boundary, and the $G$-action on $\Lambda_\mathcal{F}(G)$ is topologically free.
	\end{thm}
	
	Theorem~\ref{thm c-star simpl intro} provides further insight into the structural properties of subgroups of general type in affine buildings. We now summarize some of these properties.
	\begin{cor}
		Let $G$ be a group of general type acting properly on a locally finite affine building $X$. Then: 
		\begin{enumerate}
			\item $G$ has no amenable confined subgroups; \label{item confined}
			\item $G$ has no non-trivial uniformly recurrent subgroups. \label{item urs}
			\item $G$ satisfies Power's averaging property as defined in \cite{kennedy20}. \label{item power's}
			\item For every amenable subgroup $H $ of $G$, the quasi-regular representation $\lambda_{G/H}$ is weakly equivalent to $\lambda_G$.  \label{item weak eq}
			\item The reduced group $C^\ast$-algebra $C_r^\ast (G)$ has a unique tracial state.\label{item tracial state} 
			\item The reduced crossed product $C^\ast$-algebra $C(\Lambda_\calF(G)) \rtimes_r G$ is purely infinite simple. If $G$ is a lattice, then it is moreover nuclear. \label{item crossed product}
		\end{enumerate}
	\end{cor}
	\begin{proof}
		Properties~\ref{item confined}--\ref{item urs}--\ref{item power's}--\ref{item weak eq}--\ref{item tracial state} are well-known consequences of the $C^\ast$-simplicity of $G$ established in Theorem~\ref{thm c-star simpl intro}; see \cite{kalanthar_kennedy17,breuillard_kalantar_kennedy_ozawa17,kennedy20}.
		
		The $G$-action on $\Lambda_\calF(G)$ is minimal and topologically free by Theorem~\ref{thm G-boundary intro} and in some cases, $\Lambda_\calF(G) =\chinf$. It is moreover $n$-filling by \cite{jolissaint_robertson00}. As topological freeness implies that the action is properly outer (see \cite[Definition~1.1 \& Remark~1.3]{jolissaint_robertson00}), the reduced crossed product $C(\chinf) \rtimes_r G$ is purely infinite simple by \cite[Theorem~1.2]{robertson_steger96}. Claim~\ref{item crossed product} then follows from the fact that, in case $G$ is cocompact lattice, $\Lambda_\calF(G) = \chinf$ and the $G$-action on $\chinf$ is topologically amenable \cite[Theorem~1.1]{lecureux10}, therefore $G$ is exact and the reduced crossed product $C(\chinf) \rtimes_r G$ is nuclear \cite{anantharaman-delaroche_renault00}.
	\end{proof}
	\medskip
	
	To prove Theorem~\ref{thm G-boundary intro}, our approach relies on two main tools. The first is the dynamical behavior of strongly regular hyperbolic elements, introduced in \cite{caprace_ciobotaru15}. Using these elements, we establish the existence of a “universal” proximal sequence, which immediately implies strong proximality and is of independent interest (see Proposition~\ref{prop contraction generic position}).
	
	\begin{prop}\label{prop contraction generic position intro}
		Let \( X \) be a locally finite affine building, and let \( G \) be a group of automorphisms of \( X \) whose action on \( X \) is of general type. Then there exist strongly regular hyperbolic elements \( g_1, g_2 \in G \) such that for every ideal chamber \( C \in \Ch(X^\infty) \),
		\[
		g_2^n g_1^{n} C \xrightarrow[n \to \infty]{} C_2^+,
		\]
		where \( C_2^+ \) denotes the attracting ideal chamber of \( g_2 \).
	\end{prop}
	
	The second tool involves the study of triples of chambers at infinity. We show that there exists a subset \(\ch(\bdinf)^{[3]}_{\mathrm{gen}} \subseteq \chinf^3\) admitting an equivariant barycenter map. A triple \(\mathcal{T} \in \ch(\bdinf)^{[3]}_{\mathrm{gen}}\) is called \emph{generic}, and this terminology is justified by the fact that \(\ch(\bdinf)^{[3]}_{\mathrm{gen}}\) is open and dense in \(\chinf^3\). In the setting where \(X = \Delta^{\mathrm{BT}}(\mathbf{H}, \mathbb{K})\) is the Bruhat--Tits building of a reductive algebraic group \(\mathbf{H}\) over a non-Archimedean local field \(\mathbb{K}\), a triple of chambers is in generic position precisely when its stabilizer is a bounded subgroup of \(\bfH\).

	\begin{prop}[{Barycenter map}]
		Let \( X \) be a locally finite affine building. Then there exists a well-defined $\iso(X)$-equivariant map (called the barycenter map)
		\[\zeta: \ch(\bdinf)^{[3]}_{gen}\to X \]
		from generic triples of chambers at infinity to $X$. 
	\end{prop}
	
	Consider a probability measure $\nu $ on $\chinf$ in the measure-class as the harmonic measures $\{\nu_x\}$ defined by Parkinson in \cite{parkinson06}. These measures have full support: $\nu_x (U) >0 $ for every nonempty open set $U \subseteq \chinf$. The term ``generic'' is justified by the following result, see Proposition~\ref{prop generic full measure}. 
	\begin{prop} Let \( X \) be a locally finite affine building.
		A triple of chambers in $\chinf$ is $\nu^{\otimes 3}$-almost surely in generic position. In particular, if we let $C_1, C_2$ be opposite chambers at infinity, then the set of chambers $C_3$ such that $\{C_1,C_2,C_3\}$ is in generic position is dense in $\chinf$. 
	\end{prop}
	In \cite{parreau19}, A.~Parreau studies configurations of non-degenerate triples of chambers at infinity in $\tilde{A}_2$-buildings, and obtains a complete description of triple ratios and cross-ratios in these. One of the goals of \cite{parreau19} was to provide a geometric analogy for Fock and Goncharov coordinates \cite{fock_goncharov06}. Here we do not get such precise statements, but we provide a construction valid for affine buildings of any type and we focus on the construction and density of such generic triples of ideal chambers. It is however an interesting question to relate this barycenter map to algebraic triple ratios and cross-ratios, and to Fock and Goncharov coordinates in the context of algebraic groups. 
	
	The existence of such a barycenter map is a powerful tool for rigidity techniques, particularly when combined with the existence of $G$-equivariant measurable maps \( B \to \prob(\chinf) \), where \( B \) is a Zimmer-amenable \( G \)-space; see, for instance, \cite{bader_caprace_furman_sisto22}.  
	Our construction relies solely on the \(\mathrm{CAT}(0)\)-geometry of \( X \) and on the spherical building at infinity \( X^\infty \). In particular, straightforward adaptations show that an analogous barycenter map exists in symmetric spaces of non-compact type.
	
	\subsection{Mean proximality}
	For  an element $ w \in W_{fin}$ of the finite Weyl group $W_{fin}$ associated to $X$, the $w$-Schubert cell based at a chamber at infinity $C \in \chinf$ is 
	\[\Opp_w(C):= \{D \in \chinf \mid d_{W_{fin}}(C,D)=w\}.\]
	Let $G$ be a group of general type and let $ C \in \Lambda^+_\calF(G)$ be the attracting chamber of a strongly hyperbolic regular element.  We say that $C$ is \emph{non-avoidant with respect to Schubert cells} if we have that
	\[\Opp_w (C)\cap \Lambda_\calF(G)\neq \emptyset\]
	for all $ w \in W_{fin}$.
	
	For a probability measure $\mu $ on $G$, a measure $\nu \in \prob(\chinf)$ is $\mu$-stationary if $\mu \ast \nu = \nu $, where $\mu \ast \nu$ denotes the convolution product. We classify stationary measures associated to admissible measures on groups of general type.

	\begin{thm}[{Mean proximality}]\label{thm mean prox intro}
		Let $G <\iso(X)$ be a group of general type acting on a locally finite affine building $X$. Assume that there exists $ C \in \Lambda^+_\calF(G)$ which is non-avoidant with respect to Schubert cells. Then the $G$-action on $\Lambda_\calF(G) $ is mean proximal. In other words, for any admissible measure $\mu \in \prob(G)$, there exists a unique $\mu$-stationary measure $\nu \in \prob(\chinf)$, and this measure makes $(\chinf, \nu)$ a $(G, \mu)$-boundary. The support of $\nu$ is exactly $\Lambda_\calF(G)$. 
	\end{thm}
	
	In particular, this result applies to cocompact lattices, which have full flag limit sets, see Proposition~\ref{prop minimal}. The proof follows the ideas developed in \cite[Chapter~4]{margulis91}. In addition to establishing that \(\Lambda_\mathcal{F}(G)\) is a \(G\)-boundary, we provide in Proposition~\ref{prop equicontinuous decomposition} an equicontinuous decomposition of the flag limit set \(\Lambda_\mathcal{F}(G)\) for the \(G\)-action, which may be of independent interest.
	
	Even in the algebraic case, the property of being non-avoidant with respect to Schubert cells does not seem to be totally understood. More precisely, let $G$ be a reductive algebraic group and let $P$ be a minimal parabolic subgroup of $G$. We do not know if there exists a discrete Zariski--dense subgroup 
	$\Gamma < G$ such that its flag limit set
	$\Lambda_\calF(\Gamma) \subset G/P$ is a non-trivial proper subset of the full flag variety,
	but nevertheless intersects every Schubert cell relative to a limit flag.
	
	To the best of our knowledge, no such example is known, and no general theorem rules it out. In the case of real Lie groups, and if $\Gamma$ is Anosov, then $\Lambda_\calF(\Gamma) $ is a proper subset of $G/P$ and systematically avoids all of the Schubert cells that are not associated with the long element $w_0$. 
	In general, geometric constraints leading to a flag limit set that is not the full flag variety seems to lead to avoidance of Schubert cells.
	
	Furstenberg's boundary theory has historically been a powerful tool for studying random walks induced by group actions in the context of non-positive curvature (see, for instance, \cite{goldsheid_margulis89,benoist_quint}). In this framework, the analysis of stationary measures plays a central role. While we believe that an equicontinuous decomposition of the $G$-action on $\Lambda_\calF(G)$ holds for all groups of general type, the main application of Theorem~\ref{thm mean prox intro} in this paper is the identification of the Poisson boundaries of cocompact lattices in affine buildings, so we chose to focus on this case. We plan to address this problem for random walks on all groups of general type in a separate work.
	
	\subsection{Poisson boundaries of lattices}
	
	Finally, we identify the Poisson boundaries of cocompact lattices in affine buildings. The Poisson–Furstenberg boundary \((B_\mu, \nu_\mu)\) of a \(\mu\)-generated random walk on \(G\) is a \(G\)-space equipped with a distinguished measure class \([\nu_\mu]\), characterized as the unique maximal \((G,\mu)\)-boundary; see Theorem~\ref{thm charac poisson boundary} below or \cite{kaimanovich00} for a standard reference.
	
	The problem of finding topological models for the Poisson–Furstenberg boundary of a random walk on a group has a long and rich history. In his fundamental paper \cite{kaimanovich00}, Kaimanovich provided criteria for such identifications. We apply the strip criterion from that work to deduce the following result.
	
	\begin{thm}[{Poisson boundaries of lattices}]\label{thm poisson boundary intro}
		Let \( X \) be a thick, locally finite affine building, and let \( G \) be a discrete countable group acting properly and cocompactly on \( X \) by automorphisms. Let \( \mu \) be an admissible probability measure on \( G \), and let \( \nu \) denote the unique \(\mu\)-stationary measure on \(\chinf\) provided by Theorem~\ref{thm mean prox intro}. Then \((\chinf, \nu)\) is a compact model of the Poisson–Furstenberg boundary of \((G, \mu)\).
	\end{thm}

	\subsection{Further developments}
	The study of the reduced $C^*$-algebras of Kac--Moody groups is still largely underdeveloped compared to classical algebraic groups. We believe our dynamical and geometric approach can be helpful for better understanding the operator algebraic properties of groups acting on buildings and twin-buildings, among which Kac--Moody groups form a large class. 
	
	Moreover, $C^\ast$-simplicity can be seen as a property opposite amenability, as for amenable groups \emph{any} unitary representation is weakly contained in the regular representation. In \cite{le-bars_leibtag_vigdorovich26}, the second author shows with E.~Leibtag and I.~Vigdorovich that for some lattices in possibly exotic affine buildings of rank 2 (in particular any lattice in an affine building of type $\tilde{A}_2$), an a priori stronger instance of non-amenability is satisfied in that these group reduced $C^\ast$-algebras are selfless. Some of the arguments in \cite{le-bars_leibtag_vigdorovich26} rely on the techniques presented here.

	\begin{ackn}
		C.C is supported by a research grant (VIL53023) from VILLUM FONDEN.  C.LB is supported by ERC Advanced Grant NET 101141693, and is thankful to E.~Leibtag and I.~Vigdorovich for many helpful discussions. 
	\end{ackn}
	\section{Affine buildings}\label{section::affine_buildings}
	
	In this paper, we consider buildings as defined in \cite{abramenko_brown08}. In this paper, affine buildings $X$ are always assumed discrete (i.e. are polysimplicial complexes), with their complete system of apartments and endowed with the natural $\cat$(0) metric that is denoted by $d_X$, see \cite[Chapter~11]{abramenko_brown08}. Let us fix for the rest of this section an affine building $X$ modeled after the affine Coxeter group $W$, with associated finite Coxeter group $ W_{\mathrm{fin}}$.

	\subsection{Spherical building at infinity}
	
	As a $\cat$(0) space, $X$ has a \emph{visual} bordification, given by equivalence classes of rays, two rays being equivalent if they are at finite Hausdorff $d_X$-distance, and denoted by $\overline{X} = X \cup \bd X$. The visual boundary $\bd X$ can be endowed with a natural topology and isometries of $X$ extend to homeomorphisms on the boundary.  Moreover, when $X$ is separable, the topology on $X^\infty$ is metrizable. 
	Analogous to defining ends of trees one can equip $\bd X$ with the structure of a spherical building \cite[Propriété 1.7]{parreau00}, which we denote by $X^\infty$. We note that the map $A \mapsto A^\infty$ that associates a apartment of $X^\infty $ to an apartment $A$ of $X$ is one-to-one \cite[Proposition~4.9]{kramer_weiss14}. 
	
	\subsection{The cone topology on $\overline{X}$.} \label{section cone topo}
	
	The visual boundary is a classical notion for complete CAT(0) spaces, see e.g \cite{bridson_haefliger99}. 
	
	In an affine building \(X\), the model alcove (or chamber) may contain several special vertices; however, these vertices do not necessarily lie in the same orbit under the associated Coxeter group or weight lattice. In the literature, one such orbit is called the set of \emph{good vertices} (see \cite{parkinson06} for details). When the root system of \(X\) is non-reduced (see Section~\ref{section harmonic measures}), good vertices form a subset of the special vertices and play a key role in studying harmonic measures on \(\chinf\). If the root system is reduced, every special vertex is good. In all cases, each alcove contains at least one good vertex.
	
	After fixing a type for a special vertex, let \(X_P\) denote the set of all vertices in \(X\) of that type. For \(x \in X_P\) and a facet at infinity \(F^\infty\), there exists an affine facet based at \(x\) in the equivalence class of \(F^\infty\), denoted by \(Q(x, F^\infty)\).
	
	A detailed account on how to give a topology on the set $X \cup \chinf$, and even on $X \cup X^\infty_\tau$, where $X^\infty_\tau$ represents the set of simplices of type $\tau$ of the spherical building at infinity, is given in \cite{caprace_lecureux11}. 
	Fix a special vertex $x \in X_P$. Then the topology on $\chinf $ is given by the basis of neighborhoods $\{U_x(y) \mid y \text{ special vertex of the same type as } x \}$, where
	\begin{eqnarray}\label{eq basis chinf}
		U_x(y) := \{ C \in \ch(\bdinf) \mid  y \in Q(x, C)\} \subseteq\chinf.
	\end{eqnarray}
	
	This topology does not depend on $x$ \cite[Theorem~3.17]{parkinson06}, and makes $\chinf$ into a metrizable totally disconnected space \cite[Proposition~6.14]{grundhofer_kramer_van-maldeghem12}. In fact, the topology given by this basis is the same as the one given by the basis $\{U_x(y)\}$ \cite[Theorem~3.17]{parkinson06}, for $x,y$ special vertices of the same type. 
	
	A more detailed account on how to give a topology on boundaries of affine buildings or even masures is given in \cite{ciobotaru_muhlerr_rousseau_20}. The following proposition summarizes some properties of the compactification $X \cup \chinf$, see \cite[\S3.1]{ciobotaru_muhlerr_rousseau_20} and in particular \cite[Remark~3.4]{ciobotaru_muhlerr_rousseau_20} for the relation with the visual topology.  The set of \emph{regular points} $\bd^{\mathrm{reg}} X$ is defined as the set of boundary points $\xi \in\bd X$ that are strictly supported on a chamber, denoted $C_\xi$. 
	
	\begin{prop}
		Let $X$ be any affine building. Then there is a topology on $X \cup \chinf$ for which a basis of open sets of the chambers at infinity is given by the sets \eqref{eq basis chinf}. It agrees with the $\cat$(0) topology on $X$. This topology is first-countable and Hausdorff.
		The map $\bd^{\mathrm{reg}} X	\to \chinf$ defined by $\xi\mapsto C_\xi$ is a homeomorphism.
	\end{prop}

	Two chambers in a spherical building are said to be \emph{opposite} if the gallery distance between them is maximal, see for instance \cite[\S5.7.1]{abramenko_brown08}. In this case, there is a unique apartment joining them. Given a chamber $C_0 $ in a spherical building, the \emph{big cell} $\Opp(C_0)$ is the set of chambers opposite $C_0$. It is an open set for the cone topology.

	\subsection{Retractions}
	
	Let $X$ be an affine building, and let $C\in \chinf$. If $A$ is an apartment such that $C\subset \partial A$ (equivalently, there exists a Weyl chamber $S$ contained in $A$ representing $C$), then there exists a unique retraction map $\rho_{A, C}: X \to A$ such that $\rho_{A, C}$ preserves the distance on any apartment containing a chamber representing $C$, see \cite[Proposition 1.20]{parreau00}. Moreover, $\rho_{A, C}$ does not increase distances. We call $\rho_{A,C}$ the \emph{canonical retraction} of $X$ on $A$ based at the chamber at infinity $C \in \chinf$. In particular, for any apartment $A'$ containing $C$ in its boundary, ${\rho_{A,C}}_{\mid A'} : A' \to A$ is an isomorphism fixing $A \cap A'$ pointwise.

	\subsection{Regularity} 
	
	Let $\fraka$ be the fundamental affine apartment which we identify with the vector space on which the reflection group $W$ acts. Let $\mathfrak{a}^{+}$ be the fundamental Weyl chamber, and denote by $\mathfrak{a}^{++}$ its interior. Denote by $w_0$ the long element of the finite Weyl group $ W_{\mathrm{fin}}$ associated with the affine reflection group $W$. For $\lambda \in \mathfrak{a}$, the opposition involution $j : \mathfrak{a}^{+} \to \mathfrak{a}^{+}$ is defined by 
	$$ j (\lambda) = w_0(-\lambda).$$ 
	Let $x, y  \in X$ be distinct special vertices of the building $X$ and let $A$ be an apartment containing them. There exists a uniquely defined isomorphism sending $f: \fraka \to A $ with $f(0)= x $ and $f^{-1}(y) \in \mathfrak{a}^{+}$. The \emph{type}  (or \emph{vector distance}) between $x $ and $y$ is the vector $\theta(x,y):=f^{-1}(y) \in \fraka^+$ defined above. We say that the segment $[x,y]$ is \emph{regular} if the type $\theta(x,y)$ is regular, i.e. $\theta(x,y) \in \mathfrak{a}^{++}$.  
	
	If $[x,y]$ is regular, then both $\theta(x,y)$ and $j(\theta(y,x))$ are regular. Notice that $f(\mathfrak{a}^+)$ is a Weyl chamber of the building $X$, contained in $A$ and containing $y$. Denote by $C_y$ the chamber at infinity that it represents. Similarly, $f( \theta(x,y)- \mathfrak{a}^+)$ is a Weyl chamber containing $x$ representing a chamber $C_x \in \chinf$. These chambers are opposite in the spherical building at infinity. In this case, $A$ is the unique apartment joining them, i.e. such that both chambers belong to $A^\infty$. 
	
	\subsection{Strongly regular hyperbolic elements}\label{section SRH}
	
	We borrow the notion of strongly regular hyperbolic elements in affine buildings from \cite{caprace_ciobotaru15}. The following proposition roughly states that strongly regular hyperbolic elements satisfy a weak version of North-South dynamics, similar to loxodromic isometries in hyperbolic spaces or contracting isometries in $\cat$(0) spaces.
	
	\begin{prop}[{\cite[Proposition 2.10]{caprace_ciobotaru15}}]\label{prop dyn srh}
		Let $g \in \iso(X)$ be a type preserving strongly regular hyperbolic element, with unique translation apartment $A$. Let $C^-, C^+ \in \chinf$ be the repelling (resp. attracting) chamber at infinity for $g$. Then for every $C \in \chinf$, the limit $\lim g^n (C) $ exists and coincides with $\rho_{A, C^-} (C)$, where $\rho_{A, C^-}$ is the retraction onto $A$ centered at $C^-$. 
	\end{prop}

	\subsection{Residue building based at a vertex}
	
	Let $F$ and $F'$ be two affine facets based at a vertex $o \in X$. We say that $F$ and $F'$ have the same \emph{germ at $o$} if their intersection is an open neighborhood of $o$ in both $F$ and $F'$, see \cite[\S4.14]{kramer_weiss14}. The set of all germs at $o$ that are not reduced to $\{o\}$ can be given the structure of a simplicial spherical building, called the \emph{residue building} $\Sigma_o X$ based at $o$. It has the same type as the spherical building at infinity. Simplices of maximal dimension of $\Sigma_o X$ are the local alcoves.  
	
	For a special vertex $x \in X_P$, we have a canonical morphism of simplicial complexes 
	\begin{eqnarray}
		\Sigma_x : \bdinf \to \Sigma_x X \label{morph spherical residue}
	\end{eqnarray}
	sending any facet at infinity $F^\infty $ to $\germ_x(Q(x,F^\infty))$.

	Fix $o \in X$. For every $ x \in X$ such that $[o,x]$ is regular, there exists $C \in \chinf $ such that $x $ belongs to the interior of $ Q(o, C)$. We can then associate a unique local alcove in $\Sigma_o X$ defined by $\Sigma_o(x):=\Sigma_o (C)$. 

	In this paper, we shall consider thick affine buildings, endowed with a complete apartment system. Recall that a building is called \emph{thick} if every non-maximal simplex is contained in at least 3 distinct chambers. A special point $p \in X$ is  called thick if the spherical residue building $\Sigma_pX$ is thick. We highlight the following easy lemma. 
	
	\begin{lem}\label{lem id bd id X}
		Let $X$ be an affine building with its complete apartment system, and assume that the spherical building at infinity $X^\infty$ is thick. Let $g \in \iso(X)$ be an isometry of $X$ that fixes pointwise $\chinf$. Then $g$ acts as the identity on $X$. 
	\end{lem}
	
	\begin{proof}
		Since \( g \) acts as the identity on \( \chinf \), it acts trivially on every apartment \( A^\infty \) of the spherical building at infinity \( X^\infty \). By the one-to-one correspondence between apartments of \( X^\infty \) and those of \( X \), the restriction \( g_{|A} \) to any apartment \( A\) of \( X \) preserves \( A \) and is an automorphism of it.
		
		Because \( g \) stabilizes every apartment of \( X \), it must also stabilize thick walls -- i.e., walls that are the intersection of at least three affine apartments. As the building \( X^\infty \) is thick, there exists a thick point \( p \) (see \cite[Lemma~4.20]{kramer_weiss14}).
		
		Let \( x \in X \) be any point, and let \( A \) be an apartment containing both \( p \) and \( x \). Since \( p \) is thick, every wall in \( A \) that contains \( p \) is thick (see \cite[Lemma~4.18]{kramer_weiss14}), and hence stabilized by \( g \). As \( p \) is a vertex of \( X \), it lies at the intersection of finitely many walls, all stabilized by \( g \), which implies that \( g \) fixes \( p \).
		
		Moreover, since \( g \) stabilizes \( A \), fixes \( p \), and acts trivially on \( A^\infty \), it must fix every point of \( A \). Therefore, \( g(x) = x \). As this holds for arbitrary \( x \in X \), the isometry \( g \) must be the identity.
		
	\end{proof}
	
	If no additional assumption is specified, all the affine buildings we shall consider in this paper have thick spherical building at infinity, so that Lemma~\ref{lem id bd id X} applies. 
	
	\subsection{Harmonic measures on $\chinf$}\label{section harmonic measures}
	
	The set of chambers at infinity \(\chinf\) carries a family of measures \(\{\nu_x\}\), called \emph{harmonic measures}, indexed by the special vertices of \(X\) of a fixed type, i.e., vertices in \(X_P\) (see \cite[pp.~587–588]{parkinson06}). The interested reader can find further information on these in \cite{remy_trojan21} or \cite{bader_furman_lecureux23} (where the notion of prouniform measures is introduced). For any two special vertices \(x,y \in X\) of the same type, the measures \(\nu_x\) and \(\nu_y\) are mutually absolutely continuous \cite[Theorem~3.17]{parkinson06}. For \(\lambda \in \fraka^+\), let \(V_\lambda(x)\) denote the set of \(z \in X\) such that \(\theta(x,z) = \lambda\). Moreover 
	\[\nu_x (U_x(y)) = \frac{1}{N_\lambda} \]
	if $\theta ( x,y ) = \lambda $ and $N_\lambda= |V_\lambda(x)| $ denotes the cardinal of the set $V_\lambda(x)$ of special vertices $z\in X$ with $\theta(x,z) = \lambda$ and of the same type as $x$. In fact, we do not need such a precise information, just that the measures $\{\nu_x\}$ define a class of probability measures and that any one of them is full support: $\nu_x(U)>0$ for every open non-empty set $U$.   
	The following fact is from \cite[Proposition~2.13]{horbez_huang_lecureux23}. 
	\begin{prop}\label{prop opp full measure}
		Let $x \in X_P$, and let $C_0 \in \chinf$. Then $\nu_x$-almost every $C \in \chinf$ is opposite $C_0$. 
	\end{prop}
	\begin{cor}\label{cor opp dense}
		For every $C_0 \in \chinf$, the big cell $\Opp(C_0)$ is open dense in $\chinf$.  
	\end{cor}
	\begin{proof}
		Let $C_0\in \chinf$. That $\Opp(C_0)$ is open is classical. Assume that $\Opp(C_0)$ is not dense: there exists an open set $U \subseteq \chinf - \Opp(C_0)$. As for $x \in X_P$, the set $\{U_x(y) \mid y \in X_P \}$ forms a basis of open sets, there exists $y \in X_P$ such that $U_x(y) \subseteq \chinf - \Opp(C_0)$. But then $\nu_x(U_x(y)) >0$ and thus $\nu_x(\Opp(C_0)) <1$, contradicting Proposition~\ref{prop opp full measure}.
	\end{proof}
	The next proposition will be useful in order to construct strongly regular elements with a control on the attracting/repelling chambers. 
	\begin{prop}\label{prop restriction open sets}
		Let $U_1, U_2 \subseteq \chinf$ be open sets. Then there exists an apartment $A$ and $x,y$ special vertices in $A$ of the same type such that: 
		\begin{enumerate}
			\item $\theta(x,y) $ is regular; 
			\item $U_x(y) \subseteq U_1$ and $U_y(x) \subseteq U_2$. 
		\end{enumerate}
	\end{prop}
	
	\begin{proof}
		Let $C_1 \in U_1$. By Corollary\ref{cor opp dense}, there exists $C_2 \in U_2$ opposite $C_1$. Let $A$ be the unique apartment joining them, and let $x \in A$ be a special vertex. As $U_1, U_2$ are open, there exist special vertices $y_1,y_2$ of the same type as $x$, such that $y_1 $ is contained in the interior of the sector $ Q(x, C_1)$, resp. $y_2 $ is contained in the interior of $Q(x, C_2)$ with $U_x(y_1) \subseteq U_1$, resp. $U_x(y_2)\subseteq U_2$, and with $x$ on the geodesic segment from $y_1$ to $y_2$. With this choice, we have that $U_{y_1}(y_2) $ is contained in $U_x (y_2)$, resp. $U_{y_2}(y_1) $ is contained in $U_x (y_1)$, which proves the proposition. 
	\end{proof}
	
	As a side result, we obtain the following corollary. 
	\begin{cor}\label{cor complete opposite subsets}
		Let $U_1, U_2 \subseteq \chinf$ be open sets. Then there exist open sets $O_1 \subseteq U_1$ and  $O_2 \subseteq U_2$ such that for all $C_1 \in O_1$ and $C_2 \in O_2$, the chambers $C_1$ and $C_2$ are opposite. 
	\end{cor}
	
	\begin{proof}
		Reproduce the construction given in the proof of Proposition~\ref{prop restriction open sets}, and keep the same notations. For any $D_1 \in U_{y_2}(y_1)$ and $D_2 \in U_{y_1}(y_2)$, the alcoves $\Sigma_x (D_1)$ and $\Sigma_x(D_2)$ are opposite at $x$, hence $D_1 $ and $D_2$ are opposite. 
	\end{proof}
	
	\section{Topological and measurable boundaries}
	This section is dedicated to gathering notions from topological dynamics and boundary theory. 
	\subsection{Furstenberg boundaries}\label{section furst bd}
	
	In this section, we consider a topological dynamical system $(M, G)$, where $G$ denotes a locally compact second countable (semi)group, $M$ is a compact, often metrizable, space (the phase space) and $G \curvearrowright M$ is an action by self-homeomorphisms. 
	
	Recall that the space $\prob (M) $ of Borel probability measures on a Polish space $M$ is endowed with the (metrizable and separable) weak-$\ast$ topology.
	\begin{Def}[Minimality, proximality]
		\begin{enumerate}
			\item The space $M$ is called \emph{$G$-minimal} if every $G$-orbit is dense: for every $m \in M$, and any non-empty open subset $U \subseteq M$, there exists $g \in G$ such that $gm \in U$. 
			\item The $G$-action on $M$ is called \emph{proximal} if for any two points $m_1, m_2 \in M$, there exists a point $p \in M$, such that for every neighborhood $U$ of $p$, there exists $g \in G$ such that $gm_1, gm_2 \in U$. 
			\item The $G$-action on $M$ is called \emph{strongly proximal} if given any probability measure $\eta \in \prob (M)$, there exists $m \in M$ and
			a sequence $g_n \in G$ such that $g_n \eta \to \delta_m $. Equivalently, the $G$-action on $\prob(M)$ is proximal \cite[Lemma~III.1.1]{glasner76}. 
		\end{enumerate}
		A minimal strongly proximal $G$-space $M$ is called a \emph{$G$-boundary}. 
	\end{Def}

	A result of Furstenberg states that every group $G$ admits a $G$-boundary $\partial_F G$ which is universal in the following sense: every $G$-boundary is a continuous $G$-equivariant image of $\partial_F G$. The space $\partial_FG$ is called the \emph{Furstenberg boundary} of $G$. However the Furstenberg boundary of a discrete countable group is often ``big'': it is always extremally disconnected \cite[Proposition~2.4]{breuillard_kalantar_kennedy_ozawa17}, and therefore is not second countable unless it is reduced to a point, which is the case if and only if $G$ is amenable. 
	
	We recall that the action of a discrete countable group $G$ on a compact $G$-space $M$ is \emph{topologically free} if for all $g \in G$, the set of fixed points of $g$ 
	\[M^g := \{x \in M \mid gx= x \} \]
	has empty interior. 
	
	In this paper, we use the following breakthrough result of $C^\ast$-simplicity by Kalantar and Kennedy \cite{kalanthar_kennedy17}, see also \cite[Theorem~3.1]{breuillard_kalantar_kennedy_ozawa17}. 
	\begin{thm}\label{thm BKKO}
		Let $G$ be a discrete group. Then the following are equivalent: 
		\begin{enumerate}
			\item $G$ is $C^\ast$-simple; 
			\item $G$ acts freely on its Furstenberg boundary; 
			\item there exists a topologically free $G$-boundary. 
		\end{enumerate} 
	\end{thm}
	
	\subsection{Stationary measures}\label{section boundary theory}
	
	Let $G$ be a discrete countable group, and let $\mu$ be a probability measure on $G$.
	Let $(\Omega, \mathbb{P}) $ be the probability space $(G^{\mathbb{N}}, \mu^{\otimes \mathbb{N}})$, with the product $\sigma$-algebra. The space $\Omega$ is called the \textit{space of forward increments}. The application 
	\begin{equation*}
		(\omega=(\omega_i)_{i\in \mathbb{N}},n) \in  \Omega \times\mathbb{N} \mapsto Z_n(\omega) = \omega_1 \omega_2 \dots \omega_n,
	\end{equation*}
	defines the random walk on $G$ generated by the measure $\mu$. We say that $\mu \in \prob(G) $ is \emph{admissible} if its support $\supp(\mu)$ generates $G$ as a semigroup. 
	
	Let $Y$ be a topological space, endowed with the $\sigma$-algebra of Borel sets. Let $G \curvearrowright Y$ be a continuous action on $Y$, and let $\nu \in \prob(Y)$ be a probability measure on $Y$. We define the convolution probability measure $\mu \ast \nu$ as the image of $\mu \times \nu $ under the action $G \times Y \rightarrow  Y $. In other words, for $f $ a bounded measurable function on $Y$, 
	
	\begin{equation*}
		\int_Y f(y)d(\mu \ast \nu)(y) = \int_G\int_Y f(g \cdot y ) d\mu(g) d\nu (y).
	\end{equation*}

	\begin{Def}
		A probability measure $\nu \in \prob(Y) $ is \emph{$\mu$-stationary} if $\mu \ast \nu = \nu $. In this case, we say that $(Y, \nu)$ is a \emph{$(G, \mu)$-space}. 
	\end{Def}
	Notice that if $Y$ is a Polish space, the set of probability measures $\prob(Y)$ can be endowed with the metrizable weak-$\ast$ topology. 
	
	\begin{rem}
		When $G$ acts by homeomorphisms on a compact metric space $Y$, weak-$\ast$ compacity of $\prob(Y)$ implies that for any measure $\mu$ on $G$, there exists a $\mu$-stationary measure on $Y$. The set of $\mu$-stationary measures is convex compact in $\prob(Y)$. 
	\end{rem}
	The following statement, due to H.~Furstenberg, is one of the fundamental results of boundary theory. 
	
	\begin{thm}[{\cite{furstenberg73}}]\label{thm furstenberg73}
		Let $G$ be a locally compact second countable group,  $\mu \in \prob(G)$ an admissible measure and $(Z_n(\omega))$ be the random walk on $G$  associated to the measure $\mu$. Let $Y$ be a Polish space on which $G$ acts by homeomorphisms, and let $\nu$ be a $\mu$-stationary measure on $Y$. Then, there exists an essentially well-defined $G$-equivariant measurable map $\omega \mapsto \nu_\omega \in \prob(Y)$ such that $\mathbb{P}$-almost surely, $Z_n (\omega)\nu \rightarrow \nu_\omega $ in the weak-$\ast$ topology. Moreover, we have the decomposition $$\nu = \int_{\Omega} \nu_\omega d \mathbb{P}(\omega). $$
	\end{thm}
	
	\begin{Def}[$(G, \mu)$-boundary, mean proximality]
		\begin{enumerate}
			\item A $(G, \mu)$-space $(Y, \nu)$ is a \textit{$(G,\mu)$-boundary} if for $\mathbb{P}$-almost every $\omega\in \Omega$, the limit measure $\nu_\omega$ from Theorem~\ref{thm furstenberg73} is a Dirac measure. 
			\item A $G$-space  $Y$ is called \emph{$\mu$-proximal} if for any $\mu$-stationary measure $\nu \in \prob(Y)$, the $(G, \mu)$-space $(Y, \nu)$ is a $(G, \mu)$-boundary. 
			\item A $G$-space $Y$ is \emph{$G$-mean proximal} if for any admissible measure $\mu \in \prob(G)$, the $G$-space $Y$ is $\mu$-proximal.
		\end{enumerate}
	\end{Def}

	\begin{Def}[Equicontinuous decomposition]
		Let $(Y,d)$ be a metric $G$-space. A set $L \subseteq G$ is equicontinuous on $U \subseteq Y$ if for every $\varepsilon >0$, there exists $R >0$ such that if $y,y' \in U$ satisfy $d(y,y') < R $, then  $d(gy,gy') < \varepsilon$ for all $g \in L$. 
		
		The $G$-action on $Y$ is said to admit an \emph{equicontinuous decomposition} if we can write $G$ as a finite union of subsets $G = \cup^N_{i = 1} G_i$ and find non-empty open subsets $Y_i \subseteq Y$ such that:
		\begin{enumerate}
			\item $G Y_i = Y$, for every $i \in \{1,...,N\}$;
			\item $G_i$ is equicontinuous on $Y_i$. 
		\end{enumerate} 
	\end{Def}

	We shall use the following result as our main criterion to prove mean proximality. A statement close in spirit can be found in \cite[Proposition~VI.2.13]{margulis91}. 
	
	\begin{prop}[{\cite[Theorem~4.7]{nevo_sageev13}}] \label{thm nevo-sageev mean prox}
		If the $G$-action on $M$ is strongly proximal and admits an equicontinuous
		decomposition, then for every admissible probability measure $\mu$ on $G$, there exists a unique $\mu$-stationary measure $\nu$ on $M$, and $(M, \nu)$ is a $(G, \mu)$-boundary. In particular, $M$ is $G$-mean proximal.  
	\end{prop}
	
	The Poisson--Furstenberg boundary of a $\mu$-generated random walk on $G$ is a $G$-space with a distinguished measure class. It can be characterized by the following universality property due to Furstenberg \cite{furstenberg63}, see also \cite[Theorem~9.2]{nevo_sageev13}. 
	\begin{thm}\label{thm charac poisson boundary}
		The Poisson--Furstenberg boundary $(B_\mu , \nu_\mu)$ associated to a probability measure $\mu$ on a group $G$ is the unique maximal $(G, \mu)$-boundary: every $(G, \mu)$-boundary is a measurable $G$-factor of $(B_\mu , \nu_\mu)$.
	\end{thm} 
	Among other features, $B_\mu$ can be used to describethe set $\operatorname{Har}^\infty(G, \mu)$ of bounded $\mu$-harmonic functions on $G $  via the surjective isometry 
	\[F \in L^\infty(B_\mu, [\nu_\mu]) \mapsto \Phi_\varphi(F) \in \operatorname{Har}^\infty(G, \mu),\]
	where 
	\[\Phi_\varphi(F) (g) = \int_{B_\mu} F(gb)d\nu_\mu (b).\]
	for all $g \in G$.
	
	\section{Groups of general type}\label{section general type}
	In this section, we introduce a class of groups acting on buildings with nice proximal properties. 
	
	\subsection{Flag limit set}

	Let $G < \iso(X)$ be group acting faithfully by isometries on an affine building $(X,d)$, where $d$ is the $\cat$(0) metric. In the sequel, we will deal with groups containing strongly regular hyperbolic elements and study their dynamics on the boundary, so it makes sense to consider the flag limit set defined below. 
	\begin{Def}[Flag limit set]
		We set $\Lambda^{+}_\calF (G) \subseteq \chinf$ to be the set of attracting chambers of strongly regular hyperbolic elements in $G$. We define the flag limit set $\Lambda_\calF (G) \subseteq \chinf$ as the closure of $\Lambda^{+}_\calF (G) $ for the cone topology. 
	\end{Def}	
	
	\begin{rem}
		Classically, the limit set of a group acting $H$ on a $\cat$(0) space $Y$ is the intersection $\overline{H \cdot o} \cap \bd Y$ of the closure in $Y \cup \bd Y$ of an $H$-orbit with the visual boundary. The relation between this $\cat$(0)-limit set and ours seems subtle, and we will not use this notion. 
	\end{rem}
	\subsection{Groups of general type}
	
	Here we define the main notion of this section. 
	
	\begin{Def}\label{def general type}
		We say that a subgroup $G < \iso(X)$ is \emph{of general type} if $G$ contains strongly regular hyperbolic elements and if for any non-empty $G$-invariant closed set $M \subseteq \chinf $ and any finite set of chambers $C_1, \dots, C_N \in \chinf$, the intersection 
		\[M\cap \underset{i= 1, \dots,N}{\bigcap} \Opp(C_i)\] is non-empty. 
	\end{Def}
	
	In real algebraic groups, ``most'' (in the Baire-category sense or for the Haar measure) pairs of elements $\{g_1,g_2\}$ generate a subgroup of general type for the action on the associated symmetric space, see for instance \cite{breuillard2003dense}. We shall see in Proposition~\ref{prop minimal} that cocompact lattices in affine buildings are of general type. 
	
	We now give important properties of the flag limit set associated to groups of general type.

	\begin{thm}\label{thm flag general type}
		Let $X$ be a locally finite affine building and let $G < \iso(X)$ be a subgroup of general type. Then the limit set $\Lambda_\calF(G)$ is the smallest closed non-empty $G$-invariant subset of $\chinf$. It is $G$-minimal, perfect and uncountable. If $\Lambda_\calF(G) \neq \chinf$, then $\Lambda_\calF(G)$ has empty interior. 
	\end{thm}
	
	\begin{proof}
		First, let $C^+ $ be the attracting chamber of a strongly regular hyperbolic element $g \in G$, and let $h \in G$. Then $h g h^{-1}$ is a strongly regular hyperbolic element with attracting chamber $h C^+$. Therefore $\Lambda_\calF(G)$ is $G$-invariant. Suppose now that there exists a closed $G$-invariant subset $M\subseteq \chinf$ with
		\[\Lambda_\calF(G) - (\Lambda_\calF(G) \cap M) \neq \emptyset.\] 
		Let $g$ be a strongly regular hyperbolic element with attracting chamber $C^+$ in $\Lambda_\calF(G) - (\Lambda_\calF(G) \cap M)$, and let $C^-$ be its repelling chamber. As $G$ is of general type, $M \cap \Opp(C^-) \neq \emptyset$, and therefore there exists a chamber  $D \in M \cap \Opp(C^-)$ for which $g^n D \to C^+ $. As $M$ is closed and $G$-invariant, $C^+ \in M$, a contradiction. This shows that $\Lambda_\calF(G)$ is the smallest closed $G$-invariant subset of $\chinf$, and that it is $G$-minimal.

		The same proof shows that no point in $\Lambda_\calF(G)$ is isolated, as $\Lambda_\calF(G) \cap \Opp(C^-) \cap \Opp(C^+) \neq \emptyset$, for any strongly regular hyperbolic element $g \in G$ with attracting and repelling chambers $C^+,C^-$. As it is closed by construction, it is perfect. As it is moreover Polish and non-empty it is therefore uncountable. 
		
		Suppose finally that $\Lambda_\calF(G) \neq \chinf$, then the frontier $\partial \Lambda_\calF(G) = \overline{\Lambda_\calF(G)} -\operatorname{Int}(\Lambda_\calF(G))$ of $\Lambda_\calF(G)$ is a non-empty closed $G$-invariant set, so it must be all of $\Lambda_\calF(G)$. Therefore $\Lambda_\calF(G)$ has empty interior. 
	\end{proof}

	\section{Barycenter map}\label{section::opposite_ideal_ch}
	In this section, we construct a barycenter map for triples of ideal chambers at infinity that are pairwise opposite. We retain the notation introduced in Section~\ref{section::affine_buildings}.

	\begin{Def}\label{def generic}
		Let \( X \) be a locally finite affine building, and let \( X^\infty \) denote its ideal boundary. We say that a triple of distinct ideal chambers \( \{C_1, C_2, C_3\} \subset \Ch(X^\infty) \) is \emph{antipodal} if the chambers \( C_1, C_2, C_3 \) are pairwise opposite.
		
		We say that the triple \( \{C_1, C_2, C_3\} \) is in \emph{generic position} if it is antipodal and, for all \( 1 \leq i \neq j \leq 3 \), denoting by \( A_{ij} \) the unique apartment in $X$ having both  \( C_i \) and \( C_j \) as chambers in its boundary at infinity \(A_{ij}^{\infty} \), we have:
		\[
		A_{12}^{\infty} \cap  A_{13}^{\infty} \cap A_{23}^{\infty} = \emptyset.
		\]
		Let us denote by $\Ch(X^\infty)^{[3]}$ the set of triples of ideal chambers of $X$ that are antipodal and by $\Ch(X^\infty)_{gen}^{[3]}$ those that are in generic position.
	\end{Def}
	\begin{rem}
		Let $g \in \Aut(X^\infty)$ an automorphism of the spherical building at infinity. If $\calT \subseteq \chinf$ is an antipodal (resp. generic) triple of chambers, then $g\calT$ is still antipodal (resp, generic). 
	\end{rem}
	The following proposition relies on the convexity of the metric in $\cat$(0) spaces. 
	
	\begin{prop}
		There is a map $\zeta : \chinf^{[3]}_{gen} \to X$ which is $\iso(X)$ equivariant. 
	\end{prop}
	\begin{proof}
		Let \( \calT = \{C_1, C_2, C_3\} \subseteq \chinf \) be a triple in generic position and as always denote by \( A_{ij} \) the unique apartment in $X$ having both  \( C_i \) and \( C_j \) as chambers in its boundary at infinity \( \partial A_{ij} \). Define the function
		\begin{eqnarray}
			F_\calT : x \in X \mapsto \sum_{i\neq j} d(x, A_{ij})  \in \bbR_+\nonumber
		\end{eqnarray}

		The map $F_\calT$ is well-defined. As for $i\neq j$ the apartment $A_{ij}$ is a closed (hence complete) geodesically convex subset, the function $x \in X \mapsto d(x, A_{ij})$ is a continuous, convex function \cite[Corollary~II.2.5]{bridson_haefliger99}. As a sum of convex function, $F_\calT$ is then continuous, convex and it is bounded from below. Set $\alpha:= \inf F_\calT$ to be its infimum, and let $B(F_\calT, 1)$ be the set
		\[\{x \in X \mid F_\calT(x) \leq \alpha +1\}.\] 
		
		We claim that the set $B(F_\calT, 1)$ is bounded. 
		Assume by contradiction that it is unbounded. Up to taking a subsequence we can assume that there exists a sequence $(x_n) $ in $B(F_\calT, 1)$ that converges to a point of the boundary $\xi \in \bd X$. As $(x_n)$ is at a uniformly bounded distance from the apartments $A_{ij}$, the boundary point $\xi $ must belong to the intersection 
		$A^\infty_{12} \cap A^\infty_{13} \cap A^\infty_{23}$, which is empty as the triple is generic. This is a contradiction. 
		
		Therefore the set $B(F_\calT, 1)$ is bounded, so by local compactness of $X$ it is contained in a compact set. In particular, $F_\calT$ attains its minimum, and by convexity of $F_\calT $ this minimum is attained on a closed, geodesically convex set. Denote this minimal set by $\Min(F_\calT)$. 
		Now, in any complete $\cat$(0) space, there is a map $\operatorname{circ}$ associating the center (also called circumcenter) to any bounded subset \cite[Proposition II.2.7]{bridson_haefliger99}. As $X$ is locally finite, it is complete, and we can then consider the map $\zeta : \chinf^{[3]}_{gen} \to X$ associating to any triple $\calT \in \chinf^{[3]}_{gen} $ the center of the closed bounded set $\Min(F_\calT)$. By convexity of $\Min(F_\calT)$, and since the projection on a geodesically convex subset does not increase distances  \cite[Proposition II.2.4(4)]{bridson_haefliger99},  this center must actually belong to $\Min(F_\calT)$. This map is clearly $\iso(X)$-equivariant. 
	\end{proof}

	\begin{Def}
		The map $\zeta: \ch(\bdinf)^{[3]}_{gen}\to X $ defined above is called the \emph{barycenter map} on $\chinf$. 
	\end{Def}
	
	The following lemma makes it possible to construct triples of chambers in generic position explicitly. 
	
	\begin{lem}\label{lem construction generic triple}
		Let $C_1, C_2$ be opposite chambers at infinity. Let $A^\infty\subseteq X^\infty$ be the apartment at infinity  spanned by $C_1 $ and $C_2$. Then for any chamber 
		\[C_3 \in \underset{D \in \ch(A^\infty)}{\bigcap} \Opp(D),\] the triple $\{C_1, C_2, C_3\}$ is in generic position. 
	\end{lem}
	\begin{proof}
		Let $C_1, C_2$ be opposite chambers at infinity, and let $A_{12}^\infty\subseteq X^\infty$ be the apartment at infinity  spanned by $C_1 $ and $C_2$. Let $C_3$ be a chamber at infinity that is opposite $D$ for all $D \in \ch(A_{12}^\infty)$. As usual denote by $A_{13}^\infty$, resp. $A_{23}^\infty$, the apartment at infinity spanned by the chambers $C_1$ and $C_3$, resp. spanned by $C_2$ and $C_3$. First, we claim that the intersection $A_{12}^{\infty} \cap  A_{13}^{\infty}$ must be exactly $C_1$. Indeed, the intersection is a convex subcomplex containing $C_1$. Suppose it contains another chamber at infinity $D$. On the apartment $A_{13}^\infty$, the opposition map associating to a chamber the chamber opposite to it is a bijection \cite[Lemma~5.111]{abramenko_brown08}. Therefore the only chamber opposite $C_3$ in $A^\infty_{13}$ is $C_1$, and thus $D$ can not be opposite $C_3$. This is in contradiction with our assumption. 
		
		Of course, the same proof shows that $A_{12}^{\infty} \cap  A_{23}^{\infty} = C_2$. We then obtain immediately that $A_{12}^{\infty} \cap  A_{13}^{\infty} \cap A_{23}^{\infty} = \emptyset $, thus the triple $\{C_1, C_2, C_3\}$ is in generic position. 
	\end{proof}

	The next proposition motivates the term ``generic'' for triples of chambers in Definition \ref{def generic}. 
	
	\begin{prop}\label{prop generic full measure}
		Let $\nu$ be a probability measure in the same measure class as the harmonic measures $\{\nu_x\}$, and let $C_1,C_2 \in \chinf$ be opposite ideal chambers. Then for $\nu$-almost every $C_3 \in \chinf$, the chambers $\{C_1,C_2, C_3\}$ are in generic position. In particular, the set of $C_3$ such that the triple of chambers $\{C_1,C_2, C_3\}$ is in generic position is dense in $\chinf$. 
	\end{prop}
	\begin{proof}
		By Proposition~\ref{prop opp full measure}, for any chamber $D \in \ch(A^\infty)$, the big cell $\Opp(D)$ has $\nu$-measure 1. Therefore the intersection of finitely many big cells has measure 1, and 
		\[\nu\big(\underset{D \in \ch(A^\infty)}{\bigcap} \Opp(D)\big)= 1.\]
		By Lemma~\ref{lem construction generic triple}, any chamber $C_3 \in\bigcap_{D \in \ch(A^\infty_{12})} \Opp(D)$ is such that the triple $\{C_1, C_2, C_3\}$ is generic. This proves the first statement. 
		
		The last part of the proof is a consequence of the fact that harmonic measures have full support, as in the proof of Corollary~\ref{cor opp dense}: since every open subset $U\subseteq \chinf$ has positive $\nu$-measure, there is no such open set not containing a chamber $C_3$ for which the triple $\{C_1,C_2, C_3\}$ is in generic position. 
	\end{proof}
	Let us highlight the following immediate corollary. 
	\begin{cor}\label{cor generic triple open set}
		For every open set $U \subseteq \chinf $, there exists a triple of chambers $\{C_1, C_2, C_3\}$ in $U$ in generic position. 
	\end{cor}
	\begin{proof}
		By Corollary~\ref{cor complete opposite subsets}, there exists an opposite pair of chambers in $U$. Then the result follows from Proposition~\ref{prop generic full measure}: the set of chambers $C_3$ such that $\{C_1, C_2, C_3\} $ is in generic position is dense in $\chinf$, thus in particular $U$ contains such a chamber. 
	\end{proof}
	
	\begin{rem}
		In Lemma~\ref{lem construction generic triple} and Proposition~\ref{prop generic full measure}, we only used the $\cat$(0) geometry of the locally compact space $X$ and the spherical building at infinity $X^\infty$. In particular, the same construction works in symmetric spaces of non-compact type. 
	\end{rem}

\section{$C^\ast$-simplicity}
This section is dedicated to the proof of Theorem~\ref{thm c-star simpl intro}. 

\subsection{Cocompact lattices}
In this section, we show that cocompact lattices of affine buildings are groups of general type for which the flag limit set is actually all of $\chinf$.  Throughout this section, we assume that the groups \( G \leq \operatorname{Aut}(X) \) under consideration are not necessarily type-preserving on \( X \), and therefore not on \( X^\infty \).

The argument follows from a careful examination of \cite[Section~2]{caprace_ciobotaru15}. Specifically, cocompactness implies that there are only finitely many orbits of special vertices, which allows us to apply \cite[Proposition~2.9]{caprace_ciobotaru15}.

\begin{prop}\label{prop minimal}
	Let $G$ be a discrete countable group acting cocompactly by automorphisms on an irreducible, locally finite affine building $X$. Then  $\Lambda_\calF(G) = \chinf $ and $G$ is of general type. In particular,  the $G$-action on $\chinf$ is minimal.
\end{prop}

\begin{proof}
	By \cite[Theorem~1.2]{caprace_ciobotaru15}, the group \( G \) admits strongly regular hyperbolic elements.
	
	To prove the remainder of the corollary, it suffices to show that \( \Lambda^{+}_{\mathcal{F}}(G) \) is dense in \( \chinf \). Indeed, let \( x \in X \) be a special vertex. Choose \( C \in \chinf \) and let
	\[
	\gamma : [0,\infty] \to X
	\]
	be a strongly regular hyperbolic ray starting at \( \gamma(0) := x \) and ending at \( \gamma(\infty) := \xi \in X^{\infty} \), where \( \xi \) lies in the interior of the ideal chamber \( C \). Along the ray \( \gamma \), consider a sequence \( \{x_n\}_{n \geq 0} \) of special vertices of the same type as \( x \), with \( x_0 = x \), such that the distance \( d_X(x_n, x_{n+1}) \) is constant for all \( n \geq 0 \). Since \( G \) acts cocompactly on \( X \), the set \( \{x_n\}_{n \geq 0} \) is contained in finitely many \( G \)-orbits. Each orbit is discrete because \( X \) is a locally finite simplicial complex.
	
	Therefore, we can apply \cite[Proposition~2.9]{caprace_ciobotaru15}. There exists an increasing sequence \( \{f(n)\}_{n \geq 0} \) of positive integers such that, for any \( m > n > 0 \), there is a strongly regular hyperbolic element \( g_{n,m} \in G \) satisfying:
	\[
	x_{f(n)}, x_{f(m)} \in \Min(g_{n,m}) \quad \text{and} \quad g_{n,m}(x_{f(n)}) = x_{f(m)}.
	\]
	Denote the attracting ideal chamber of \( g_{n,m} \) by \( C_{n,m}^+ \). Observe that one can construct a strongly regular geodesic ray
	\[
	\gamma_{n,m} : [0,\infty] \to X
	\]
	starting at \( \gamma_{n,m}(0) = x \), passing through the geodesic segment \( [x_{f(n)}, x_{f(m)}] \), and ending at an ideal point \( \xi_{n,m} \) in the interior of \( C_{n,m}^+ \). Moreover, the geodesic segment \( [x, x_{f(m)}] \) is contained in \( \gamma \cap \gamma_{n,m} \), and the ray \( [x_{f(n)}, \xi_{n,m}) \) lies in \( \Min(g_{n,m}) \).
	
	Taking \( n,m \to \infty \), we obtain that the sequence \( \{C_{n,m}^+\}_{m > n} \subset \Lambda^{+}_{\mathcal{F}}(G) \) converges to \( C \) with respect to the cone topology on \( \chinf \). The result follows.
\end{proof}

As a direct corollary of Proposition~\ref{prop minimal}, we obtained the following useful result, which we highlight separately for future reference. 
\begin{cor}\label{cor srh prescribed open sets}
	Let $X$ be an irreducible, locally finite affine building, and let $G$ be a group of automorphisms of $X$  acting cocompactly on $X$. Let $U\subseteq \chinf $ be an open set. Then there exists a strongly regular hyperbolic isometry $g \in G$ with attracting chamber in $U $. 
\end{cor}

\subsection{Strong proximality}

In this section, we prove the existence of a universal proximal sequence for groups of general type. 

\begin{prop}\label{prop contraction generic position}
	Let \( X \) be a locally finite affine building. Consider two strongly regular hyperbolic elements \( g_1, g_2 \in \Aut(X) \). For each \( g_i \), denote its attracting and repelling chambers by \( C_i^+ \) and \( C_i^- \), respectively. Let \( A_1 \) be the unique affine apartment in \( X \) determined by the pair \( (C_1^+, C_1^-) \), and let \( A_1^\infty \) denote its ideal boundary. Assume that 
	\[
	C_2^- \in \bigcap_{D \in \ch(A_1^\infty)} \Opp(D).
	\]
	Then for every ideal chamber $C \in \chinf$, we have
	\[g_2^{n} g_1^{n} C  \xrightarrow[n \to \infty]{} C_2^+.\]
\end{prop}
\begin{proof}
	Let $A_2$ be the unique affine apartment in $X$ spanned by $C^+_2$ and $C^-_2$. As well, since $C_2^-$ is opposite any $D \in \ch(A_1^{\infty})$, let $A_D$ be the unique affine apartment in $X$ spanned by $C^-_2$ and $D$. Then, fix a point $x_2$ in $A_2$ such that $x_2$ lies in the interior of a common sector corresponding to the ideal chamber $C_2^-$ of the apartments $A_2$ and $A_D$, for each $D \in \ch(A_1^{\infty})$.
	
	Again by the hypothesis that  $C_2^-$ is opposite any $D \in \ch(A_1^{\infty})$, and by \cite[Proposition 2.10]{caprace_ciobotaru15}, we have that 
	\[
	g_2^{n}(D) \xrightarrow[n \to \infty]{} C_2^+,
	\]  
	for any $D \in \ch(A_1^{\infty})$.
	
	For what follows, fix $r_2 >0$ with $r_2 \in \mathbb{R}$, the geodesic ray \( \rho_2 := [x_2, \xi_2^+)  \subset A_2\), where \( \xi_2^+ \) denotes the barycenter of the ideal chamber \( C_2^+ \), and the standard open neighborhood \( U(x_2, r_2, C_2^+) \) of \( C_2^+ \), centered at the base point \( x_2 \), with radius \( r_2 \), and associated to the geodesic ray \( \rho_2 \). Then there is 
	$$N_{r_2}>0 \text{ such that } g_2^{n}(D)  \in U(x_2, r_2, C_2^+),$$ 
	for any $n \geq N_{r_2}$ and any  $D \in \ch(A_1^{\infty})$.
	
	Let now $x_1$ be a point in $A_1$. For each $D \in \ch(A_1^{\infty})$ consider the geodesic rays  \( \rho_{D,1} := [x_1, \xi_D)  \subset A_1\) and  \( \rho_{D,2} := [x_2, \xi_D)  \subset A_D\), where \( \xi_D \) denotes the barycenter of the ideal chamber \( D \). Since \( X \) is a \( \mathrm{CAT}(0) \) space, the geodesic rays $\rho_{D,1}, \rho_{D,2}$ are in fact parallel. So, for each  $D \in \ch(A_1^{\infty})$, we fix a point $y_D \in \rho_{D,1} \cap A_D$, as far as we need on the geodesic $\rho_{D,1}$, and such that $y_D$ is in the interior of the common sector corresponding to the ideal chamber $D$ of the apartments $A_1$ and $A_D$.
	
	\medskip
	Let $C \in \chinf$. By  \cite[Proposition 2.10]{caprace_ciobotaru15}, there exists a unique $D' \in \ch(A_1^\infty)$ such that $g_1^m (C) \xrightarrow[m \to \infty]{}  D'$. Let \( U(x_1, y_{D'}, D') \) be the standard open neighborhood of \(D'\), centered at the base point \( x_1 \), with radius \( y_{D'} \), and associated to the geodesic ray \(  \rho_{D',1} \). Then there is 
	$$M_{C,1}>0 \text{ such that } g_1^m (C) \in  U(x_1, y_{D'}, D'),$$ 
	for every $m \geq M_{C,1}$. Notice that the constant $M_{C,1}$ depends on the ideal chamber $C$. 
	
	Since we can choose $y_{D'}$ as far as we need on the geodesic $\rho_{D,1}$, and independently on the ideal chamber $C$, we have that a large part of the geodesic segment $[x_1, y_{D'}]$, say the geodesic subsegment $[z_{D'}, y_{D'}]$, is contained in the interior of the common sector corresponding to the ideal chamber $D'$ of the apartments $A_1$ and $A_{D'}$. 
	
	Then, notice that for any $n \geq N_{r_2}$, we have 
	$$C_{2}^- \in \ch(g_2^n(A_{D'})), \; \; x_2, g_2^{n}(x_2) \in g_2^n(A_{D'}) \text{ and } g_2^{n}([z_{D'}, y_{D'}]) \subset g_2^n(A_{D'}).$$
	
	Patching geodesics by using the geodesic rays 
	$$[x_2, g_2^n(\xi_{D'})), \quad [g_2^n(x_2), g_2^n(\xi_{D'})),$$ 
	$$[g_2^n(z_{D'}),g_2^n(\xi_{D'})),  \text{ and }  [g_2^n(z_{D'}),  g_2^ng_1^{m}(\xi_{C})),$$
	one constructs the geodesic ray $[x_2, g_2^ng_1^{m}(\xi_{C}) )$ and notices that 
	$$ g_2^ng_1^{m}(C) \in U(x_2, r_2, C_2^+), \text{ for all } m \geq M_{C,1}, n\geq N_{r_2}.$$
	
	Taking any $n \geq \max\{M_{C,1}, N_{r_2}\}$, we obtain that $g_2^ng_1^{n}(C) \in U(x_2, r_2, C_2^+)$, and the conclusion of the proposition follows.
\end{proof}


\begin{prop}\label{prop contraction general type}
	Let $G <\iso(X)$ be a subgroup of general type. Then there exists a sequence $\{g_n\}_{n \geq 0} $ in $G$ and a chamber $C^+\in \Lambda_\calF(G)$ such that for any $C \in \chinf$, 
	\[g_n C  \xrightarrow[n \to \infty]{} C^+.\]
\end{prop}

\begin{proof}
	Let $g_1$ be a strongly regular hyperbolic element with attracting chamber $C^+_1$ and repelling chamber $C^-_1$. Let $A^\infty$ be the apartment at infinity spanned by $C^+_1$ and $C^-_1$. As $G$ is of general type and $\Lambda_\calF(G)$ is closed $G$-invariant by Theorem~\ref{thm flag general type}, the intersection
	\[ \Lambda_\calF(G) \bigcap_{D \in \ch(A^\infty)} \Opp(D)\]
	is non-empty. As the set $\bigcap_{D \in \ch(A^\infty)} \Opp(D)$ is open, there exists a strongly regular hyperbolic element $g_2 \in G$ with repelling chamber $C_2^- $ in it. Then by Proposition~\ref{prop contraction generic position}, for any chamber $C \in \chinf$, 
	\[g_2^{n} g_1^{n} C  \xrightarrow[n \to \infty]{} C_2^+,\]
	where $C_2^+$ is the attracting chamber of $g_2$. This proves the claim. 
\end{proof}

In the proof of Proposition~\ref{prop contraction general type}, we obtained the existence of a generic triple in $\Lambda_\calF(G)$. 
\begin{cor}\label{cor exist generic triple general type}
	Let $G < \iso(X)$ be a subgroup of general type. Then there exists a triple of chambers $\calT \subseteq \Lambda_\calF(G)$ in generic position.
\end{cor}
\begin{proof}
	Reproduce the same construction of the proof of Proposition~\ref{prop contraction general type}. Then the triple $\{C_1^-, C_1^+, C_2^-\}$ is in generic position by Lemma~\ref{lem construction generic triple}. 
\end{proof}

We get the following as an immediate corollary. 
\begin{cor}\label{cor strong prox general type}
	Let $G <\iso(X)$ be a subgroup of general type. Then the $G$-action on $\chinf $ is strongly proximal. 
\end{cor}
\begin{proof}
	Let $\nu \in \prob(\chinf)$ be a probability measure. Consider the universal proximal sequence $\{g_n\}_{n \in \mathbb{N}} \subset G$ given by Proposition~\ref{prop contraction general type}. Let $f \in \mathcal{C}(\chinf)$ be a continuous function. Then
	\[ g_n \nu(f) = \int f(g_n C) d\nu(C). \]
	As for all chamber $C \in \chinf$, the sequence $g_n C $ converges to the chamber $C^+$, we have by the dominated convergence theorem applied to the sequence of functions $\{f_n(\cdot):=f(g_n \cdot)\}$ on $\chinf$ that 
	\[g_n \nu(f) \to f(C^+).\]
	Therefore, $g_n \nu$ converges to the Dirac measure $\delta_{C^+}$ in the weak-$\ast$-topology, as $n \to \infty$. 
\end{proof}

We can now refine Corollary~\ref{cor exist generic triple general type}. 

\begin{cor}\label{cor everywhere generic triple flag limit set}
	Let $G <\iso(X)$ be a subgroup of general type, and let $U \subseteq \chinf $ be an open set intersecting $\Lambda_\calF(G)$. Then there exists a generic triple of chambers at infinity $\calT$ entirely contained in $\Lambda_\calF(G) \cap U$. 
\end{cor}

\begin{proof}
	There exists a generic triple $\calT':= \{C_1, C_2, C_3\}$ of chambers at infinity contained in $\Lambda_\calF(G)$ by Corollary~\ref{cor exist generic triple general type}. By strong proximality of the $G$-action on $\chinf$, Proposition \ref{prop contraction generic position}, there exists a sequence $\{g_n\}_{n\geq 0}$ and a chamber $D \in \Lambda_\calF(G)$ such that $g_n \calT' \to D$. As the $G$-action on $\Lambda_\calF(G)$ is minimal, $D$ can be chosen in $U$, proving the claim. 
\end{proof}
\subsection{Topological freeness}
Recall that the action of a discrete group $H$ on a compact space $Y$ is topologically free if the set of fixed points $Y^h$ of any non-trivial element $h \in H$ has empty interior.

\begin{prop}\label{prop top free general type}
	Let $G < \iso (X)$ be a subgroup of general type and assume that the action $ G \curvearrowright X$ is proper. Let $h \in G$. If there exists an open set $U\subseteq \chinf $ such that $\Lambda_\calF(G) \cap U \neq \emptyset$ with $h_{|U\cap\Lambda_\calF(G)}= \Id$, then $h = \id$. In other words, the $G$-action on $\Lambda_\calF(G)$ is topologically free.
\end{prop}

\begin{proof}
	Assume that there exists an open set $U \subseteq \chinf^h $ intersecting $\Lambda_\calF(G)$ with $h_{|U\cap\Lambda_\calF(G)}= \Id$. There exists a generic triple $ \calT \in \chinf^3$ in $U \cap \Lambda_\calF(G)$ by Corollary~\ref{cor everywhere generic triple flag limit set}. If we denote by $z:= \zeta(\calT)$ its barycenter, we have that $h (z) = z$. As $\Lambda_\calF^+(G)$ is dense in $\Lambda_\calF(G)$ by definition, there exists a strongly regular hyperbolic element $g \in G$ with attracting chamber $C^+$ in $U$, and repelling chamber $C^-\in \Lambda_\calF(G)$. In particular, for all $n $ large enough, we can assume that $g^{-n}h g^n \zeta(\calT) = \zeta(\calT) = z$. As the action is proper, there exists $a \in G$ such that $g^{-n}h g^n = a$ for infinitely many $n $. If $a = \Id $, we obtain that $h = \id$. If not, Proposition~\ref{lem id bd id X} implies that there exists a non-empty open set $V \subseteq \chinf - \chinf^a$. Now, as $\Opp(C^-) $ is dense in $\chinf $, the set $\Opp(C^- ) \cap V$ is non-empty, and for any chamber $D \in\Opp(C^- ) \cap V$ we have that $g^{-n} h g^n  D = D$ for all $n $ large enough. In particular $aD= D$, contradicting the fact that $D \in \chinf - \chinf^a$.
\end{proof}
\subsection{$C^\ast$-simplicity}
\begin{proof}[Proof of Theorem~\ref{thm c-star simpl intro} and Theorem~\ref{thm G-boundary intro}]
	The set $\Lambda_\calF(G)$ is non-empty compact as it is a closed subset of the compact set $\chinf$. The $G$-action on $\Lambda_\calF(G)$ is minimal by Theorem~\ref{thm flag general type}  and strongly proximal by Corollary~\ref{cor strong prox general type}. The flag limit set $\Lambda_\calF(G)$ is then a $G$-boundary. The action on $\Lambda_\calF(G)$  is moreover  topologically free by Proposition~\ref{prop top free general type}. This proves Theorem~\ref{thm G-boundary intro}.  We can then apply Theorem~\ref{thm BKKO}, finishing the proof of Theorem~\ref{thm c-star simpl intro}. 
\end{proof}

\section{Equicontinuous decomposition and mean proximality}

In this section, we prove Theorem~\ref{thm mean prox intro}. We let $X$ be an irreducible locally finite affine building, and let $G \curvearrowright X$ be an action of a group $G$ by automorphisms. 

\subsection{Uniformities and their bases}

We recall a few facts from general topology regarding uniformities on sets. Let $Y$ be a set. Its diagonal $\Delta(Y)$ is the set $\{(y,y) \mid y \in Y\}$. An \emph{entourage} on $Y$ is a set $V \subseteq Y\times Y $ containing $\Delta(Y)$ and such that $V = V^{-1}$, where $Y^{-1} := \{(y_2,y_1) \mid (y_1,y_2) \in V\}$. We denote the set of entourages on $Y$ by $\calD_Y$. 

For a set $V \subseteq Y\times Y$, we define the following subset of $Y\times Y$
\[V \circ V:=\{(x,z) \mid \exists \ y \in Y \ : \ (x,y) \in V \text{ and } (y,z) \in V\}.\]
\begin{Def}
	A \emph{uniformity} $\mathcal{U}$ on $Y$ is a family of entourages of $Y$ such that 
	\begin{enumerate}[label=(U\arabic*)]
		\item If $V \in \mathcal{U}$ and $V \subseteq W \in \calD_Y$, then $W \in \calU$. 
		\item If $V_1$, $V_2 \in \calU$, then $V_1 \cap V_2 \in \calU$. 
		\item For every $V \in \calU$, there is $W \in \calU$ such that $W\circ W \subseteq V$. (This is called the “uniform triangle inequality”.)
		\item $\cap_{V \in \calU} V = \Delta(Y)$. 
	\end{enumerate}
	Otherwise saying, a uniformity on $Y$ is a filter on the set $Y \times Y$ whose elements are entourages.
	
	A \emph{basis}  for the uniformity $\calU$ is a family  $\calB \subseteq \calU$ such that for any $ V\in \calU$, there exists $W \in \calB$ with $W \subseteq V$.
\end{Def}
We shall construct uniformities using the following result. 

\begin{thm}[{\cite{engelking68}}]\label{thm engelking uniformity}
	Let $\calB \subseteq \calD_Y$ be a family of entourages on $Y$ satisfying the following conditions:
	\begin{enumerate}[label=(BU\arabic*)]
		\item If $V_1$, $V_2 \in \calB$, there exists $W \in \calB$ with $W \subseteq V_1 \cap V_2$. 
		\item For every $V \in \calB$, there is $W \in \calB$ such that $W \circ W	 \subseteq V$. 
		\item $\cap_{V \in \calB} V = \Delta(Y)$. 
	\end{enumerate}
	Then the family $\calU \subseteq \calD_Y$  of entourages containing an element of $\calB$ is a uniformity, for which $\calB$ is a basis. 
\end{thm}
For an element $V \in \calU$ of a uniformity $\calU$ on $Y$, and $x \in Y$, the ball $B(x,V)  $ denotes the set of all $y \in Y$ with $(x,y ) \in V$. Given that, a uniformity $\calU$ on a space $Y$ induces a topology on $Y$ that is defined by the collection of  “open sets”:
\[\{ O \subseteq Y \mid \forall x \in O, \ \exists \ V \in \calU \ : \ B(x, V) \subseteq O \}.\]

\begin{ex}[Metric spaces]
	Conversely, when $(Y,d)$ is a metric space, the sets $V_i := \{(x,y) \mid d(x,y) < 2^{-i}\} $ define a basis for a uniformity on $Y$. The topology induced by this uniformity is the same as the one induced by metric $d$. 
\end{ex}

\subsection{Equicontinuous decomposition}

The goal of this subsection is to prove the following proposition. 
\begin{prop}\label{prop equicontinuous decomposition}
	Let $G < \iso (X)$ be a subgroup of general type. Assume that $\Lambda_\calF(G)$ contains an apartment at infinity. Then the $G$-action on $\Lambda_\calF(G)$ admits an equicontinuous decomposition. 
\end{prop}

Let us introduce some notations. Let \(\fraka\) be a fundamental apartment of \(X\), regarded as the model apartment of \(X\). Let \(\fraka^+\) be a fundamental Weyl chamber in \(\fraka\), regarded as the model Weyl chamber of \(X\).

For what follows, fix a special vertex $o\in X$. Recall that as in Section~\ref{section cone topo}, a basis of open sets in $\ch (\bdinf)$ is given by the sets $U_{o}(y)$ for which  $\theta (o, y)$ is colinear with the barycenter vector $\rho$ in the model Weyl chamber \(\fraka^+\) of $X$. We adopt this basis from now on. In particular, when we write $U_{o}(y)$, for $o,y $ special vertices, we always mean that $\theta (o, y)$ is colinear with the barycenter vector $\rho$ in \(\fraka^+\). 
Let $\widetilde{d^\infty}$ be a metric on $\chinf$ compatible with the topology induced by the open sets $U_{o}(y)$. 

\begin{rem}
	The construction of the metric $\widetilde{d^\infty}$ is given for instance in the unpublished lecture notes of F.~Paulin \cite[\S3.7, Proposition 3.70]{paulin_notes}. One can have in mind that the function \(\widetilde{d^\infty}\) is reminiscent of visual metrics in Gromov-hyperbolic spaces.
\end{rem}

For \(r \ge 0\), define \(\mathcal{U}(r)\) to be the collection of open sets of the form \(U_o(y)\), where, according to the convention above, the direction \(\theta(o,y)\) is colinear with the barycenter vector $\rho$ in \(\fraka^+\), and \(d_X(o,y) \ge r\).

\begin{lem}\label{lem uniformity topology}
	For every $\varepsilon >0$, there exists $r \geq 0$ such that if $C, D \in U_o(y)$, with $U_o(y) \in \calU(r)$, then $\widetilde{d^\infty}(C, D) \leq \varepsilon$. 
\end{lem}
\begin{proof}
	Let us assume the contrary. Consider sequences $U_n:= U_{o, y_n} \in \calU(n)$ and $C_n , D_n \in U_n $ such that $\widetilde{d^\infty}(C_n, D_n) \geq \varepsilon$ for all $n$. By compacity of $\chinf$, up to taking a subsequence we can assume that there exists $C_\ast \in \chinf$ such that $C_n \to_n C_\ast $ in the cone topology. Let $U= U_{o} (y_\ast)$ be an open neighborhood of $C_\ast$. By convergence, we can choose $n_0\geq d(o, y_\ast)$ large enough such that $C_n \in U_{o, y_\ast}$ for all $n \geq n_0$. In particular, for $n \geq n_0$, both $y_\ast $ and $y_n$ belong to $Q(o, C_n)$. Since we chose $y_\ast $ and $y_n$ to be such that $\theta(o, y)$ and $\theta(o, y_n)$ are barycentric in $Q(o, C_n)$, the geodesic segment $[o, y_n]$ contains $[o, y_\ast]$ for all $n\geq n_0$. Since $y_n \in Q(o, D_n)$ for all $n$, it implies that for $n \geq n_0$, we also have $y_\ast \in Q(o,D_n)$. Since this is true for all neighborhood $U_{o, y_\ast}$ of $C_\ast$, we get that $D_n$ converges to $ C_\ast$ as well. This contradicts the fact that $\widetilde{d^\infty}(C_n, D_n) \geq \varepsilon $ for all $n$ as the metric $\widetilde{d^\infty}$ is compatible. 
\end{proof}

\begin{figure}[h]
	\centering
	\begin{center}
		\begin{tikzpicture}[scale=0.4]
			\draw (-8, 4) -- (8,-4)  ;
			\draw (0, 4) -- (16,-4)  ;
			\draw (-8, 2) -- (4,-4)  ;
			\draw (-8, -2) -- (4,4)  ;
			\draw (-8, 0) -- (0,-4)  ;
			\draw (-8, 2) -- (-4,4)  ;
			\draw (-8, -2) -- (-4,-4)  ;
			\draw (-8, 0) -- (0,4)  ;
			\draw (0, -4) -- (14,3)  ;
			\draw (-4, 4) -- (12,-4)  ;
			\draw (4, 4) -- (16,-2)  ;
			\draw (4, -4) -- (16,2)  ;
			\draw (8, -4) -- (16,0)  ;
			\draw (8, 4) -- (16,0)  ;
			\draw (-4, -4) -- (12,4)  ;
			\draw (12, 4) -- (16,2)  ;
			\draw (12, -4) -- (16,-2)  ;
			\draw (-8, -4) -- (8,4)  ;
			\draw (-8, 0) -- (16,0)  ;
			\draw (-8, 3) -- (14,3)  ;
			\draw (-8, 1) -- (16,1)  ;
			\draw (-8, 2) -- (16,2)  ;
			\draw (-8, -1) -- (16,-1)  ;
			\draw (-8, -3) -- (16,-3)  ;
			\draw (-8, -2) -- (16,-2)  ;
			\draw (0, 4) -- (1,6)  ;
			\draw (16, -4) -- (19,-1)  ;
			\draw (4, 2) -- (5.2,6)  ;
			\draw (8, 0) -- (19,2)  ;
			\filldraw (0,0) circle(5pt) ;
			\filldraw (-6,-1) circle(5pt) ;
			\filldraw (8,2) circle(5pt) ;
			\draw (0,0) node[above]{$o$} ;
			\draw (-6,-1) node[below left]{$y$} ;
			\draw (8,2) node[above right]{\(g^{-1}o\)} ;
			\draw (-2,-0.5) node{\(c_y\)} ;
			\draw [thick, fill=yellow, opacity=0.2]
			(0, 0) -- (8,4) -- (12,4)  -- (14,3) -- (16,2) -- (16,0);	
			\draw [thick, fill=yellow, opacity=0.2]
			(4, 2) -- (5.2,6) -- (19,2)  -- (8,0);	
			\draw [thick, fill=gray, opacity=0.2]
			(0,0) -- (-4,0) -- (-2,-1)  ;	
		\end{tikzpicture}
	\end{center}
	\caption{The alcove \(c_y\) and a subset of $\overline{\Opp_o(c_y)}$ in yellow. }\label{figure sets G(y)}
\end{figure}
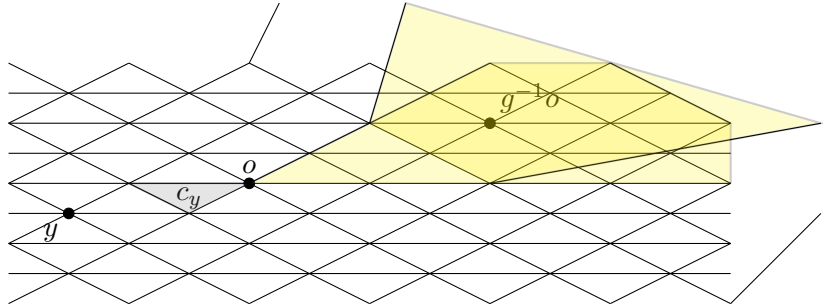

For $r \geq 0$, define $\calD(r) \subseteq \chinf\times \chinf$ to be the set of pairs of chambers $(C_1, C_2)$ such that there exists $U \in \calU(r) $ with $C_1, C_2 \in U$. As a union of open sets, $\calD(r)$ is an open set in $\chinf\times \chinf$, with respect to the product topology. It clearly contains the diagonal. 
\begin{cor}\label{cor uniformity}
	The family of sets $\calD(r)$, with $r \geq 0$, define a basis for a uniformity on $\chinf$ which is compatible with the cone topology on $\chinf$. In particular, any open neighborhood $U$ of the diagonal with respect the product topology is contained in some $\calD(r)$. 
\end{cor}
\begin{proof}
	The family $\calD(r)$ satisfies the conditions of Theorem~\ref{thm engelking uniformity}. In particular, it defines a basis for a uniformity on $\chinf$. By Lemma~\ref{lem uniformity topology}, this uniformity is compatible with the metric $\dinf$ whose induced topology corresponds to the cone topology on $\chinf$.
\end{proof}
For an alcove $c \in \Sigma_o X$, let $\overline{\Opp_o(c)}$ be the union of the closed sectors $\overline{Q(o, C)}$ with base point $o$ such that $\Sigma_o (C)$ and $c$ are opposite in $\Sigma_o X$. Let $y \in X$ be such that $\theta(o, y)$ is regular and let us denote by $c_y:= \Sigma_o (y)$ the alcove in the residue of $o$ corresponding to the projection of $y$ on the residue building $\Sigma_o X$. Let
\[\mathcal{G}(y):= \{g \in G \mid g^{-1}o \in \overline{\Opp_o(c_y)} \} \subseteq G .\]
In Figure~\ref{figure sets G(y)}, we represent a situation where $g \in \mathcal{ G }(y)$. 

\begin{lem}\label{lem equicontinuous act}
	Let $y \neq o $ be a special vertex such that $\theta(o, y)$ is colinear with the barycenter vector $\rho$ in \(\fraka^+\). Then $\mathcal{G}(y)$ acts equicontinuously on $U_o(y)$. 
\end{lem}

\begin{proof}
	Let us write $r := d(o, y)>0$.  Let $g \in \mathcal{G}(y)$, and let $C, D \in U_o(y)$. As $y$ and $g^{-1} o $ belong to the interior of opposite Weyl chambers based at $o$, we have that $y \in Q(g^{-1} o, C)$ and $y \in Q(g^{-1} o,D)$. 
	Therefore, we have that $gy \in Q(o, gC) $ and $gy \in Q(o, gD) $. Moreover, $d_X(o, gy) = d_X(g^{-1}o, y) \geq d_X(o, y)$ since $\angle_o(y, g^{-1}o) \geq \pi/2$ and by the law of cosines (observe that for any two points $w,z $ in a Weyl sector based at $o$, we have that $\angle_o(w, z) \leq \pi/2$). As a consequence, $gC$ and $gD $ both belong to $U_o(y) \in \calD(r)$. The proof then follows from Corollary \ref{cor uniformity}.
\end{proof}

The following result is standard but useful. 
\begin{lem}\label{lem ch opp apartment}
	Let $A$ be an apartment in a spherical building $\Delta$. Then every chamber $C\in \ch(\Delta)$ is opposite a chamber in $A$. 
\end{lem}
\begin{proof}
	Let $C$ be a chamber in $\Delta$. By \cite[Theorem~(A.19)]{ronan89}, the building $\Delta$ is homotopic to a bouquet of spheres, the number of spheres being equal to the number of chambers opposite $C$. Consider the geometric realization of $\Delta$. By the proof of \cite[Theorem~(A.19)]{ronan89}, the building $\Delta'$ obtained by deleting in $\Delta$ all the chambers opposite $C$ is contractible. In particular, it can not contain the apartment $A$ as the latter is homotopic to a sphere. 
\end{proof}
For a special vertex $o \in X$ and an apartment $A$ containing $o$, we denote by $\Sigma_oA$ the corresponding apartment in the residue building $\Sigma_oX$ and by $\ch(\Sigma_oA)$ its chambers. These chambers are in one-to-one correspondence with the alcoves in $A$ containing $o$. To ease the notations, we make this identification. Here as in the rest of the paper, we make the assumption that $X$ is thick, and thus $o$ can be taken to be a thick point.
\begin{lem}\label{lem partition G}
	Consider a thick special vertex $o \in X$ and let $A$ be an apartment containing $o$. For every alcove $c \in \ch(\Sigma_oA)$, we let $y \in X$ be such that 
	\begin{enumerate}
		\item $y$ belongs to the interior of $c$; 
		\item $\theta(o,y)$ is colinear with $\rho$. 
	\end{enumerate}
	Let $S$ be the finite collection of such $y$. Then we have that \[\bigcup_{y \in S} \mathcal{ G } (y)= G.\]
\end{lem}
\begin{proof}
	We need to prove that
	\[\bigcup_{c \in \ch(\Sigma_oA)}  \overline{\Opp_o(c)} = X.\]
	By Lemma~\ref{lem ch opp apartment}, any chamber $c' \in \ch(\Sigma_oX)$ is opposite a chamber in $\Sigma_oA$. As a consequence, for all $g \in G$, there exists $c \in \ch(\Sigma_oA)$ for which the element $g^{-1}o $ belongs to the closure of a Weyl chamber $Q(o, C)$ such that $\Sigma_o C$ and $c$ are opposite alcoves. 
	In particular, for any $g \in G$, there exists $y \in S$ with $g \in \mathcal{G} (y)$ and the result follows. 
\end{proof}

\begin{lem}\label{lem covering for cocompact}
	Assume moreover that the $G$-action on $X$ is of general type. Then for any $y \in X$ such that $U_o(y) \cap\Lambda_\calF(G)\neq \emptyset$, we have that $\Lambda_\calF(G) \subseteq G \cdot U_{o}(y)$. 
\end{lem}
\begin{proof}
	The $G$-action on $\Lambda_\calF(G)$ is minimal by Theorem~\ref{thm flag general type}. In particular, for all $U\subseteq \chinf$ non-empty open set in $\Lambda_\calF(G)$ and any $C \in \Lambda_\calF(G)$, there exists $g \in G$ with $gC \in U$. Observe that the set $U_o(y) \cap\Lambda_\calF(G)$ is open in $\Lambda_\calF(G)$ and non-empty by assumption, so that for any chamber $C \in \chinf$, there exists $g \in G$ with $gC \in U_o(y)$. The result follows. 
\end{proof}
\begin{proof}[Proof of Proposition \ref{prop equicontinuous decomposition}]
	By assumption, $\Lambda_{\mathcal F}(G)$ contains an apartment at infinite $A^\infty$. Denote by $A$ the corresponding affine apartment, and let $o $ be a special vertex in $A$. Consider as in Lemma~\ref{lem partition G} the finite set $S$ obtained by taking points $y$ belonging to interior of alcoves of $A$ in the residue of $o$. By Lemma~\ref{lem partition G}, we have that
	\[\bigcup_{y \in S} \mathcal{ G } (y)= G.\]
	By the choice of $A$, for any $y \in S$, there exists a chamber $C \in \ch(A^\infty)$ with $C \in U_o(y)$. In particular, for every $y \in S$, we have that $U_o(y) \cap \Lambda_\calF(G) \neq \emptyset$. By Lemma~\ref{lem covering for cocompact}, we have that $\Lambda_\calF(G) \subseteq G \cdot U_{o}(y)$. 
	By Lemma~\ref{lem equicontinuous act}, $\mathcal{ G }$ acts equicontinuously on $U_o(y)$. This proves the result. 
\end{proof}
\begin{cor}
	Let $X$ be a thick and locally finite affine building and let $G$ be a discrete countable group acting on $X$ properly and cocompactly by automorphisms. Then the $G$ action on $\chinf$ admits an equicontinuous decomposition. 
\end{cor}
\begin{proof}
	By Proposition~\ref{prop minimal}, $G$ is of general type and $\Lambda_{\mathcal F}(G) = \chinf$. The assumptions of Proposition~\ref{prop equicontinuous decomposition} are thus satisfied. 
\end{proof}

The result below shows that the assumption that there exists an apartment at infinity in $\Lambda_\calF(G)$ is equivalent to the existence of a chamber $C \in \Lambda^+_{\mathcal F}(G) $ which is non-avoidant with respect to Schubert cells. 
\begin{lem}\label{lem non-avoidant}
	Let $G$ be a group of general type in an affine building $X$. Let $g \in G$ be a strongly regular hyperbolic element and let $C^\pm$ its attracting and repelling chambers at infinity. Denote by $A^\pm$ its unique fixed apartment at infinity. If for some $w \in W_{fin}$, the flag limit set $\Lambda_{\mathcal F}(G)$ intersects the Schubert cell 
	\[\Opp_w(C^-) = \{ C\in \chinf \mid d_{W_{fin}}(C^-, C) = w\},\]
	then $\Lambda_\calF(G)$ contains the unique chamber $C_w \in A^\pm$ with $d_{W_{fin}}(C^-, C_w) = w$. In particular, if for every $w \in  W_{fin}$, 
	\[\Lambda_{\mathcal F}(G) \cap \Opp_w(C^-)\neq \emptyset, \]
	then $\Lambda_\calF(G)$ contains $A^\pm$. 
\end{lem}
\begin{proof}
	Let $C \in \Lambda_{\mathcal F}(G) \cap \Opp_w(C^-)$. Then $g^n C \to C_w$ as $n \to \infty$. As $\Lambda_{\mathcal F}(G)$ is closed and $G$-invariant, $C_w$ must belong to $\Lambda_{\mathcal F}(G)$.
\end{proof}

\subsection{Mean proximality}

\begin{proof}[Proof of Theorem~\ref{thm mean prox intro}]
	Lemma~\ref{lem non-avoidant} ensures that if there exists a chamber $ C \in \Lambda^+_\calF(G)$ which is non-avoidant with respect to Schubert cells, then there there exists an apartment at infinity all of whose chambers belong to $\Lambda_\calF(G)$.  By the criterion given in Proposition~\ref{thm nevo-sageev mean prox}, the result then follows from Corollary~\ref{cor strong prox general type} providing strong proximality and Proposition~\ref{prop equicontinuous decomposition}, which asserts that $G \curvearrowright \Lambda_\calF(G)$ admits an equicontinuous decomposition. 
\end{proof}

\section{Poisson boundaries of lattices}

Let $\mu$  be an admissible probability measure on a discrete countable group $G$. We say that $\mu$ is \emph{of finite first logarithmic moment} if there exists a distance function $|\cdot | : G \to \bbR_+$ on $G$ which is quasi-isometric to a word metric for which
\[\sum_{g \in G}\log |g| \mu(g)<\infty.\]
Note that being of finite logarithmic moment is independent of the word metric chosen. 
The \emph{entropy} of $\mu$ is defined by 
\[H(\mu) = \underset{n  \to \infty}{\lim} \frac{1}{n} \sum_{g \in G} \mu^{\ast n}(g) \log \mu^{\ast n}(g). \]
We denote by $\mui$ the probability measure on $G$ defined by $\mui(g) = \mu(g^{-1})$: it is still an admissible measure.

Let now $X$ be a thick and locally finite affine building and let $G$ be a discrete countable group acting on $X$ cocompactly by automorphisms. Consider on $G$ a distance function induced by this proper and cocompact action. Let $\mu$ be an admissible probability measure on $G$. By Theorem~\ref{thm mean prox intro}, there exists a unique $\mu$-stationary measure on the set of chambers at infinity that we denote by $\nu$. Similarly, there exists a unique $\mui$-stationary measure on the set of chambers at infinity that we denote by $\nui$. We denote by $(B, \nu_B)$ the $(G, \mu)$-Poisson--Furstenberg boundary and by $(\Bi, \nu_{\Bi})$ the $(G, \mui)$-Poisson--Furstenberg boundary, as introduced in Section~\ref{section boundary theory}. 

In this section, we show that if $\mu$ has finite first logarithmic moment, then $(\chinf, \nu)$ is a compact metrizable model for the Poisson--Furstenberg boundary of $(G, \mu)$. This will follow from the strip criterion due to Kaimanovich, which we recall now. 
\begin{thm}[{\cite[Theorem~6.5]{kaimanovich00}}]\label{thm kaiman strip}
	Let $\mu$ be a probability measure of finite first logarithmic moment and finite entropy on a group $G$. Let $(B, \nu)$ be a $(G, \mu)$-boundary and $(\check{B}, \nui)$ be a $(G, \mui)$-boundary. Assume that there exists a measurable $G$-equivariant map  
	\[ S : (b, \check{b}) \in B\times \Bi \mapsto S(b, \check{b}) \subseteq G.\]
	If for $\nu \times \nui$-almost all pair $ (b, \check{b})$, the strip $S (b, \check{b})$ has polynomial growth, then $(B, \nu)$ and $(\Bi, \nui)$ realize the Poisson--Furstenberg boundaries of $(G, \mu)$ and $(G, \mui)$ respectively. 
\end{thm}
The main result of this section is Theorem~\ref{thm poisson boundary intro}, which we recall now for the convenience of the reader. 
\begin{thm}\label{thm poisson boundary}
	Let $X$ be a thick and locally finite affine building and let $G$ be a discrete countable group acting on $X$ properly and cocompactly by automorphisms. Let $\mu$ be an admissible probability measure on $G$ and let $\nu$ be the unique $\mu$-stationary measure on $\chinf$ given by Theorem~\ref{thm mean prox intro}. Then $(\chinf, \nu)$ is a compact model of the Poisson--Furstenberg boundary of $(G, \mu)$. 
\end{thm}
We define 
\[\calO := \{(C_1, C_2) \in\chinf^2\mid C_1 \text{ opposite } C_2\}\] 
to be the set of pairs of opposite chambers at infinity. By Corollary~\ref{cor opp dense} and Proposition~\ref{prop restriction open sets}, it is a non-empty open subset of $ \chinf^2$ (for the product topology), and it is $\Aut(X)$-invariant. 

We begin by proving that $\calO$ has full $\nu\otimes \nui$-measure. 
\begin{lem}
	A pair of chambers $(C_1, C_2) \in \chinf^2$ is $\nu\otimes \nui$-almost surely antipodal. 
\end{lem}
\begin{proof}
	First, it is a known fact that the support of $\nu$ is a closed $G$-invariant set, and so by minimality of the $G$-action on $\chinf$ given by Proposition~\ref{prop minimal}, $\nu$ has full support, i.e. any non-empty open set $O \subseteq \chinf$ has positive $\nu$-measure. The same holds for the measure $\nui$. The $G$-action on $(B \times \Bi, \nu_{B} \otimes\nui_{\Bi})$ is ergodic by a classical theorem of Kaimanovich \cite{kaimanovich03}. By Theorem~\ref{thm mean prox intro}, the space $(\chinf, \nu) $ is a $(G, \mu)$-boundary, so by Theorem~\ref{thm charac poisson boundary} it is a measurable $G$-factor of $(B, \nu_B)$, resp. $(\chinf, \nui) $ is a measurable $G$-factor of $(\Bi, \nu_{\Bi})$. Therefore $(\chinf^2, \nu\otimes \nu) $ is a measurable $G$-factor of $(B\times \Bi, \nu_B \otimes \nu_{\Bi})$. In particular, by an easy argument, this implies that  the $G$-action on $(\chinf^2, \nu\otimes \nui) $ is ergodic as well. As $\nu$ and $\nui$ have full support, the space $\calO$ has positive $\nu\otimes \nui$-measure. As it is $G$-invariant, it has full measure by ergodicity. 
\end{proof}
We define the map $\calA\colon (C_1, C_2) \in \calO \mapsto \calA (C_1, C_2) \subseteq X$ that associates to a pair of opposite chambers the unique apartment in $X$ that they span. 
\begin{proof}[Proof of Theorem~\ref{thm poisson boundary}]
	By cocompactness of the $G$-action on $X$, consider a finite family  $(x_1, \dots, x_p)$ of vertices that are representatives of the $G$-orbits in the set of vertices of $X$. We say that a vertex $x \in X$ is of type $t_j$ if it belongs to the $G$-orbit $G\cdot x_j$. Consider the sets $\calO_j$ defined by 
	\[\calO_j := \{(C_1, C_2) \in \calO \mid \calA(C_1, C_2) \text{ contains a vertex of type  }t_j\}.\]
	It is easily seen that any $\calO_j$ is open and therefore $\nu\otimes \nui$-measurable. Indeed fix a pair of opposite chambers at infinity $(C_1,C_2) \in \chinf^2$ such that $\calA(C_1,C_2)$ contains a vertex $o$ of type $t_j$. Let $y_1, y_2 \neq o$ be special vertices on the geodesic rays $[o, \xi_{C_1}) $ and $[o, \xi_{C_2}) $ respectively, where $\xi_{C_i}$ denotes de barycenter of $C_i$ for $i= 1,2$.  Then any pair $(D_1, D_2) \in U_o(y_1) \times U_o(y_2) $ is in $\calO_j$. 
	
	We have that $\calO = \bigcup_{j=1}^p \calO_j$, so there exists $j \in \{1, \dots, p\}$ so that $\calO_j$ has positive $\nu\otimes \nui$-measure.  Without loss of generality, assume that $j= 1$. As the set $\calO_j$ is $G$-invariant, we obtain by ergodicity of the $G$-action on $(\chinf^2, \nu\otimes \nu) $ that it must have full $\nu\otimes \nui$-measure. 
	
	For any $(C_1, C_2) \in \calO_1$, define the strip  $S' (C_1,C_2)$ to be the union of all the cosets $g \operatorname{Stab}_{G}(x_1)$, with $g$ such that $g x_1 \in \calA(C_1,C_2)$. As the number of vertices in $\calA(C_1,C_2)$ has polynomial growth (it is isometric to $\bbR^d $ for some $d$), it has polynomial growth with respect to the word metric on $G$ induced by the cocompact action. Now for any such $g$, the set $g\operatorname{Stab}_{G}(x_1) \subseteq G$ is finite by properness of the action, so that $S'(C_1, C_2)$ has polynomial growth. The proof then follows from Theorem~\ref{thm kaiman strip}. 
\end{proof}

\appendix
\section{A triple of chambers in generic position}

Let $k$ be a field and let $E= (e_1,e_2,e_3)$ be the canonical basis associated to the vector space $k^3$. Let $G = \SL_3(k)$ and let $P, P' < G$ be the minimal parabolic subgroups corresponding to the full flags $\{0\} \subset \langle e_1\rangle \subset \langle e_1,e_2\rangle \subset k^3$ and $\{0\} \subset \langle e_3\rangle \subset \langle e_2,e_3\rangle \subset k^3$ respectively. In coordinates, these correspond to the subgroups of $G$ defined by

\[ P = \{\begin{pmatrix}
	a & b & c  \\
	0 & e & f  \\
	0 & 0 & i   
\end{pmatrix} \mid aei = 1\}, \text{ and }
P' = \{\begin{pmatrix}
	a & 0 & 0  \\
	d & e & 0  \\
	g & h & i   
\end{pmatrix}  \mid aei = 1\}.
\]
Notice that as $P \cap P' $ is a maximal split torus $T:= \{\operatorname{Diag}(a,e, (ae)^{-1}) \mid a,e \in k \}$, so that $P$ and $P'$ are transverse and correspond to opposite chambers in the Tits building of $G$. 

Let $v = e_1 + e_2 + e_3$ and for $t \in k$, consider the linear plane defined by 
\[V_t := \{ - xt + (1+t) y - z = 0 \mid x,y,z \in k\}.\]

Observe that for all $t$, $v \in V_t$. Consider the parabolic subgroup $P^t := \stab\{\{0\} \subset \langle v \rangle \subset V_t \subset k^3\}$. 
\begin{prop}
	For all $t \neq 0, -1$, the parabolic subgroups $P, P'$ and $P^t$ are in generic position. 
\end{prop}

We make explicit computations for $t= 1$, where $V_1 = \{ (2y - z, y, z ) \mid y, z \in k\}$. 

Then one easily computes that the vector $(2, 1, 0 )$ also belongs to $V_1$, so that if a matrix 
\[A := \begin{pmatrix}
	a & b & c  \\
	d & e & f  \\
	g & h & i   
\end{pmatrix} \]
belongs to $P^1$, then 
\[\left\{ 
\begin{array}{ll}
	a+ b +c &= d+e+f = g+h+i \\
	2a + b & = 2(2d+e) - (2g+h).
\end{array} 
\right.
\]
In particular, the intersection $P \cap P^1$ is the group
\[T_1 := \{A_{a,e} := \begin{pmatrix}
	a & 2e - 2a & (ae)^{-1} -2e +a  \\
	0 & e & (ae)^{-1}- e  \\
	0 & 0 & (ae)^{-1}   
\end{pmatrix} \mid a, e \in k \}.\]
\begin{rem}
	Observe that $(a,e) \mapsto A_{a,e}$ is a morphism from the free abelian group $k^2$ to $T_1$, so that $T_1$ is a maximal split torus.  
\end{rem}

Observe that the maximal split torus corresponding to the intersection $P \cap P'$ stabilizes the apartment corresponding to the complex given by Figure~\ref{figure apartment A_2}:

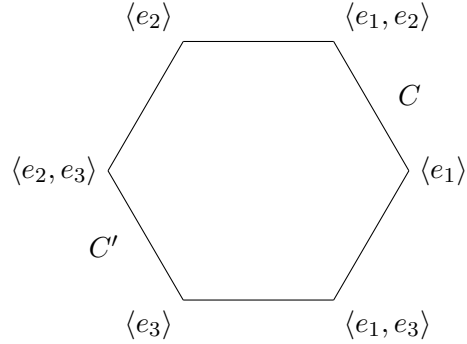
\begin{figure}[h]
	\centering
	\begin{center}
		\begin{tikzpicture}[scale=2]
			\draw (-0.5, 0.86) -- (0.5, 0.86)  ;	
			\draw (1, 0) -- (0.5, -0.86)  ;
			\draw (1, 0) -- (0.5, 0.86)  ;
			\draw (0.5, -0.86) -- (-0.5, -0.86)  ;	
			\draw (-1, 0) -- (-0.5, -0.86)  ;
			\draw (-1, 0) -- (-0.5, 0.86)  ;		
			\draw (-0.5, 0.86) node[above left]{$\langle e_2 \rangle$} ;
			\draw (0.5, 0.86) node[above right]{$\langle e_1, e_2 \rangle$} ;
			\draw (-0.5, -0.86) node[below left]{$\langle e_3 \rangle$} ;
			\draw (0.5, -0.86) node[below right]{$\langle e_1, e_3 \rangle$} ;
			\draw (1,0) node[right]{$\langle e_1 \rangle$} ;
			\draw (-1,0) node[left]{$\langle e_2, e_3 \rangle$} ;
			\draw (0.86,0.5) node[right]{$C$} ;
			\draw (-0.86,-0.5) node[left]{$C'$} ;
		\end{tikzpicture}
	\end{center}
	\caption{The apartment corresponding to $T$ in the spherical building of $G$}\label{figure apartment A_2}
\end{figure}
We highlighted the chambers $C$ and $C'$ stabilized by $P$ and $P'$ respectively. Notice that the torus $T_1$ corresponding to the intersection $P \cap P^1$ does not stabilize the partial flag $\{0\} \subset \langle e_2 \rangle \subset k^3$ nor the partial flag $\{0\} \subset \langle e_1, e_3 \rangle \subset k^3$, so that the intersection of the apartments associated to $T$ and $T_1$ is exactly the chamber $C$. 

By repeating similar computations, one can show that $P'$ and $P^1$ are also transverse, and that their intersection contains the maximal split torus 
\[T_1 := \{B_{a,e} := \begin{pmatrix}
	a & 0 & 0 \\
	a-e & e & 0  \\
	a - 2e + (ae)^{-1} & 2e - 2(ae)^{-1} & (ae)^{-1}   
\end{pmatrix} \mid a, e \in k \}.\]

Among the simplices in the apartment given by Figure~\ref{figure apartment A_2}, this torus only stabilizes the chamber associated to the flag $\{0\} \subset \langle e_2 \rangle \subset\langle e_2,e_3 \rangle \subset k^3$. 

Observe that for $t= 0 $, $P^t$ stabilizes $e_1$ while for $t= -1$, $P_t$ stabilizes $e_2$. The general argument for $t \neq 0,-1$ amounts to showing that $V_t $ does not contain these lines. In this particular case, we get a concrete uncountable family of triples in generic position. The existence of such a family is given by Proposition~\ref{prop generic full measure}.
\begin{cor}
	There exists an uncountable family of minimal parabolic subgroup $P^t$ such that $(P, P', T^t)$ is in generic position. Otherwise the vertices associated to the parabolic subgroups corresponding to the partial flags $\{0\} \subset \langle e_2 \rangle \subset k^3$ or  $\{0\} \subset \langle e_1, e_3 \rangle \subset k^3$  would be stabilized. 
\end{cor}

\bibliographystyle{alpha}
\bibliography{bibliography}
\end{document}